\documentclass[12pt]{amsart}
\usepackage{amsmath,amssymb,txfonts}
\usepackage{amsthm, color,mathrsfs}
\IfFileExists{citeref.sty}{\usepackage{citeref}}{}
\usepackage{graphicx}
\usepackage[colorlinks=true,
            linkcolor=blue,
            citecolor=blue,
            urlcolor=blue]{hyperref}

\makeatletter
\@namedef{subjclassname@2020}{%
  \textup{2020} Mathematics Subject Classification}
\makeatother
\numberwithin{equation}{section}
      \newtheorem{thm}{Theorem}[section]
      \newtheorem{definition}{Definition}

      \newtheorem{rem}[thm]{Remark}
      \newtheorem{cor}[thm]{Corollary}
      
      \newtheorem{lem}[thm]{Lemma}
      \newtheorem*{carlesoncorona}{Carleson's corona theorem}
      \newtheorem*{wolffideal}{Wolff's ideal theorem}
      \newtheorem*{treilcounterexample}{Treil's counterexample}
      \newtheorem*{slicecorona}{Corona theorem for slice regular functions}

      \newcommand{\ct}[1]{\langle {#1}\rangle \lower.3ex\hbox{$_{t}$}}
      \newcommand{\lt}[1]{[ {#1}] \lower.3ex\hbox{$_{t}$}}

\newcommand{\bd}{\mathbb{D}}

\newcommand{\rd}{\mathrm{d}}

\begin{document}

\title[Wolff-type ideal theorems for slice regular functions]
{Wolff-type ideal theorems for slice regular functions}

\author[Z. Xu]{Zhenghua Xu}
\address{Zhenghua Xu, School of Mathematics,\newline
Hefei University of Technology, Hefei, Anhui 230601, P. R. China}
\email{zhxu@hfut.edu.cn}

\author[W. Yang]{ Wenwan Yang}
\address{Wenwan Yang, School of Mathematics and Statistics,\newline
Guangdong University of Technology, Guangzhou, Guangdong 510520, P. R. China}
\email{wenwan\_yang@163.com}

\author[C. Yuan]{  Cheng Yuan$^\ast$  }
\address{Cheng Yuan, School of Mathematics and Statistics,\newline
Guangdong University of Technology, Guangzhou, Guangdong 510520, P. R. China}
\email{yuancheng1984@163.com}

\begin{abstract}
We study Wolff-type ideal problems in $H^\infty(\mathbb B)$, the algebra of
bounded slice regular functions on the quaternionic unit ball.
We construct $f,f_1,f_2\in H^\infty(\mathbb B)$ satisfying
$|f|^2\le |f_1|^2+|f_2|^2$ such that
$f^{\ast N}\notin\mathcal I(f_1,f_2)$ for every integer $N\ge1$,
where $\mathcal I(f_1,f_2)$ denotes the right ideal generated by $f_1$ and $f_2$.
We also construct intrinsic functions $a,a_1,a_2\in H^\infty(\mathbb B)$ satisfying
$|a|\le |a_1|+|a_2|\le1$ and
$a^{\ast2}\notin\mathcal I(a_1,a_2)$.
For intrinsic data satisfying $|f|\le\sum_{\ell=1}^n|f_\ell|$,
the cube $f^{\ast3}$ has a representation in the generated right ideal
with intrinsic coefficients.
For general data, we obtain sufficient conditions for ideal membership
under a fourth-order pointwise bound with a Dini factor.
If the generators preserve a common slice, a quadratic bound with
a Dini factor suffices.
\end{abstract}


\thanks{$^\ast$ Corresponding author.  }


\date{}



\maketitle


\section{Introduction}\label{s1}

Let $\mathbb D$ be the unit disk in $\mathbb C$. The bounded holomorphic functions on $\mathbb D$ form the Banach algebra $H^\infty(\mathbb D)$
under pointwise multiplication, with the supremum norm $\|f\|_\infty=\sup_{z\in\mathbb D}|f(z)|$.
 Given finitely many generators  $f_1,\ldots,f_n\in H^\infty(\mathbb D)$, the ideal
$\sum_{\ell=1}^n f_\ell H^\infty(\mathbb D)$ consists of all sums
$\sum_{\ell=1}^n f_\ell g_\ell$ with $g_\ell\in H^\infty(\mathbb D)$.

 The ideal equals the whole algebra if and only if it contains the constant function $1$.
This requires a uniform positive lower bound for
$\sum_{\ell=1}^n\left|f_\ell\right|$. Carleson's corona theorem
\cite{C}, proved in 1962, shows that this condition is also sufficient.
\begin{carlesoncorona}
Let $n\ge1$ and $f_1,\ldots,f_n\in H^\infty(\mathbb D)$.
If
\begin{equation}\label{eta1}
 \inf_{z\in\mathbb D}\sum_{\ell=1}^n\left|f_\ell(z)\right|=\delta>0,
\end{equation}
then there exist $g_1,\ldots,g_n\in H^\infty(\mathbb D)$
such that
\[
 \sum_{\ell=1}^n f_\ell g_\ell=1.
\]
\end{carlesoncorona}

A necessary condition for $g\in H^\infty(\mathbb D)$ to belong to this
ideal is the bound $|g|\le C\sum_{\ell=1}^n|f_\ell|$ for some $C>0$.
Wolff's theorem asserts that every $g\in H^\infty(\mathbb D)$ satisfying
this bound has its cube in the ideal. Rescaling $g$ gives the following
formulation.

\begin{wolffideal}
Let $n\ge1$ and $f_1,\ldots,f_n,g\in H^\infty(\mathbb D)$.
If
\begin{equation}\label{etb1}
 |g(z)|\le\sum_{\ell=1}^n\left|f_\ell(z)\right|,
 \qquad z\in\mathbb D,
\end{equation}
then there exist $g_1,\ldots,g_n\in H^\infty(\mathbb D)$
such that
\[
 g^3=\sum_{\ell=1}^n f_\ell g_\ell.
\]
\end{wolffideal}
See, for example, \cite[Chapter~VIII]{G} for these two theorems and their proofs.

 The exponent $3$ in Wolff's ideal theorem works for every finite family of
generators. Treil \cite{Tr} proved that it is the smallest possible
positive integer: the square may fail even for two generators.
\begin{treilcounterexample}
There exist $f,f_1,f_2\in H^\infty(\mathbb D)$ such that
\[
 |f(z)|\le \left|f_1(z)\right|+\left|f_2(z)\right|,
 \qquad z\in\mathbb D,
\]
but
\[
 f^2\notin f_1H^\infty(\mathbb D)+f_2H^\infty(\mathbb D).
\]
\end{treilcounterexample}

In 2006, Gentili and Struppa introduced slice regular functions of one
quaternionic variable as an extension of holomorphic functions
of one complex variable \cite{Gentili-Struppa-06,Gentili-Struppa-07}.
The theory includes convergent power series of the form
\[f(q)=\sum_{\nu=0}^{\infty}q^\nu a_\nu, \qquad q,a_\nu\in\mathbb H,\]
where  $\mathbb H$ denotes   the quaternion algebra.

The pointwise product of two slice regular functions need not be
slice regular. We therefore use the slice product (or $\ast$-product)
defined in Section~\ref{s2}. Denote by $H^\infty(\mathbb B)$ the real
Banach algebra of bounded (left) slice regular functions on the
 quaternionic unit ball $\mathbb B=\{q\in\mathbb H:|q|<1\}$, equipped with the $\ast$-product
and the supremum norm
\[ \|f\|_\infty=\sup_{q\in\mathbb B}|f(q)|.\]

In 2016, Gentili, Sarfatti and Struppa \cite{GSaS} proved that finitely
many slice regular functions with no common zeros on an axially
symmetric slice domain generate the entire ring of slice regular
functions, without boundedness assumptions on the generators or
coefficients. For bounded functions, Colombo, Pozzi, Sabadini and
Wick \cite{CPSW} proved the following analogue of Carleson's corona theorem.
\begin{slicecorona}
Let $n\ge1$, $0<\delta<1$, and
$f_1,\ldots,f_n\in H^\infty(\mathbb B)$. Suppose that
\begin{equation}\label{etd}
 \delta^2\le\sum_{\ell=1}^n\left|f_\ell(q)\right|^2\le1,
 \qquad q\in\mathbb B.
\end{equation}
Then there exist $g_1,\ldots,g_n\in H^\infty(\mathbb B)$
such that
\[
\sum_{\ell=1}^n f_\ell\ast g_\ell=1.
\]
\end{slicecorona}
The upper bound in \eqref{etd} is obtained by rescaling the generators.

In hypercomplex analysis, related matrix corona and symmetry methods appear in the work of Lin, Lu and Zu \cite{Lin} on quaternionic Hardy spaces, and in the work of  Colombo, Pozzi, Sabadini and Wick \cite{CPSWmatrix} on slice monogenic functions  with values in Clifford algebras.

Given generators $f_1,\ldots,f_n\in H^\infty(\mathbb B)$, denote the right ideal
generated by $f_1,\ldots,f_n$ by
\[
 \mathcal I\left(f_1,\ldots,f_n\right)
 =\left\{\sum_{\ell=1}^n f_{\ell}\ast g_{\ell}:
               g_1,\ldots,g_n\in H^\infty(\mathbb B)\right\},
\]
where no norm closure is taken.

In the setting of slice regular functions, our first main result shows that pointwise domination does not ensure that any positive slice power belongs to the generated right ideal.
\begin{thm}\label{Wolff-no}
\label{t1}
There exist $f,f_1,f_2\in H^\infty(\mathbb B)$ such that
\[
 |f(q)|^2\le \left|f_1(q)\right|^2+\left|f_2(q)\right|^2,
 \qquad q\in\mathbb B,
\]
but
\[
 f^{\ast N}\notin\mathcal I\left(f_1,f_2\right),
 \qquad N\in\mathbb N,
\]
where $\mathbb N=\{1,2,\ldots\}$ and $f^{\ast N}=\underbrace{f\ast f\ast \cdots \ast f}_N$.
\end{thm}
The generators in this construction preserve the slice
$\mathbb C_i=\mathbb R+i\mathbb R$, while
$f(\mathbb B\cap\mathbb C_i)\not\subseteq\mathbb C_i$.
The proof uses the matrix representation $\chi$ in \eqref{echidef}:
the upper right entry of $\chi(f)^N$ has simple zeros at a fixed sequence
for every $N\ge1$. Differentiating a putative ideal relation at these
zeros contradicts the Schwarz--Pick estimate.
Two generators are necessary for this phenomenon: for a single
generator, $|f|\le C\left|f_1\right|$ implies $f\in\mathcal I\left(f_1\right)$, as shown
at the end of Section~\ref{s2}.

We next consider intrinsic data. The \emph{intrinsic functions} form
the closed commutative real subalgebra
\[
H_+^\infty(\mathbb B)=\left\{ f \in H^\infty(\mathbb B): f\left(\bar q\right)=\overline{f(q)}, q\in \mathbb{B}\right\}
\]
of $H^\infty(\mathbb B)$.

The square can still fail to belong to the generated right ideal when all the functions are intrinsic. The following theorem is obtained from a real-symmetric version of Treil's counterexample.

\begin{thm}\label{thm2-1}
There exist $a,a_1,a_2\in H_+^\infty(\mathbb B)$ such that
\begin{equation}\label{eq:intrinsic-square-data}
 |a(q)|\le \left|a_1(q)\right|+\left|a_2(q)\right|\le1,
 \qquad q\in\mathbb B,
\end{equation}
but
\begin{equation}\label{eq:intrinsic-square-obstruction}
 a^{\ast2}\notin\mathcal I\left(a_1,a_2\right).
\end{equation}
\end{thm}

To prove Theorem~\ref{thm2-1}, Lemma~\ref{lloctreil} adapts the finite interpolation construction in \cite[Theorem~1.2 and Section~3]{Tr}.
We then reflect the data to obtain real symmetry.
Uniform lower bounds outside a boundary neighborhood preserve the failure of ideal membership under this reflection.

For intrinsic functions, Wolff's ideal theorem, combined with real symmetrization, gives the
following cubic conclusion. Theorem~\ref{thm2-1} shows that the exponent cannot be reduced to two.
\begin{cor}
\label{t2}
Let $n\ge1$ be an integer.
\begin{enumerate}
\renewcommand{\labelenumi}{\textup{(\roman{enumi})}}
\item
Let $f,f_1,\ldots,f_n\in H_+^\infty(\mathbb B)$ satisfy
\begin{equation}\label{et21}
 |f(q)|\le\sum_{\ell=1}^n\left|f_\ell(q)\right|,
 \qquad q\in\mathbb B.
\end{equation}
Then there exist $g_1,\ldots,g_n\in H_+^\infty(\mathbb B)$
such that
\begin{equation}\label{et22}
 f^{\ast3}=\sum_{\ell=1}^n f_\ell\ast g_\ell.
\end{equation}
\item
For $n\ge2$, there exist $f,f_1,\ldots,f_n\in H_+^\infty(\mathbb B)$
such that
\[
 |f(q)|\le\sum_{\ell=1}^n\left|f_\ell(q)\right|\le1,
 \qquad q\in\mathbb B,
\]
but
\begin{equation}\label{et23}
 f^{\ast2}\notin\mathcal I\left(f_1,\ldots,f_n\right).
\end{equation}
\end{enumerate}
\end{cor}

\medskip

Theorem~\ref{t1} rules out a Wolff power conclusion for general
slice regular data under pointwise domination. We therefore seek stronger
bounds that imply membership of $f$ itself in the right ideal.
For related scalar ideal problems, see
\cite{Bou,Ceg90,Ceg94,Pau05,Rao,Sa,Tr,Tr07,YY23}.
Given  $f, f_1,\ldots,f_n\in H^\infty(\mathbb B)$, denote
\[
 S(q)=\sum_{\ell=1}^n \left|f_{\ell}(q)\right|,\qquad q\in\mathbb B. \]
\medskip
\noindent\textbf{Question:  }
For which nondecreasing functions $h:[0,\infty)\to[0,\infty)$,
with $h(0)=0$, does the condition $|f(q)|\le h(S(q))$ on
$\mathbb B$ imply $f\in\mathcal I\left(f_1,\ldots,f_n\right)$?

Corollary~\ref{cor:quadratic-gauge-obstruction} shows that, for $n\ge2$,
this implication fails whenever $h(t)\ge ct^2$ for some $c>0$ and all
sufficiently small $t>0$, even for intrinsic functions.

 We call the prescribed tuple $(f,f_1,\ldots,f_n)$ the \emph{data}; it is
\emph{intrinsic} when all its functions are intrinsic.

We say that the generators $\left\{f_1,\ldots,f_n\right\}$ \emph{preserve a common slice} if there is an imaginary unit $I\in \mathbb S$ such that
\[ f_{\ell}\left(\mathbb B\cap\mathbb C_I\right)\subseteq\mathbb C_I, \qquad 1\le \ell\le n,\]
where $\mathbb S=\{q\in\mathbb H:q^2=-1\}$ and $\mathbb C_I=\mathbb R+I\mathbb R=\{x+Iy:x,y\in\mathbb R\}$. Intrinsic generators preserve every slice.

We next give sufficient conditions for ideal membership under a fourth-order
bound with a Dini factor. If the generators preserve a common slice,
a quadratic bound with a Dini factor suffices.
\begin{thm}\label{thm2}
Let $n\ge1$ be an integer.
\begin{enumerate}
\renewcommand{\labelenumi}{\textup{(\roman{enumi})}}
\item\label{gauge:general}
Let $f,f_1,\ldots,f_n\in H^\infty(\mathbb B)$, and put
\[
 S(q)=\sum_{\ell=1}^n \left|f_\ell(q)\right|,\qquad q\in\mathbb B.
\]
Let $h:[0,\infty)\to[0,\infty)$ be nondecreasing, with $h(0)=0$,
and suppose that
\[
 |f(q)|\le h(S(q)),\qquad q\in\mathbb B.
\]
Assume that there exist constants $C>0$ and $t_0\in(0,1]$ and a bounded
nonnegative nondecreasing function $\omega:(0,1]\to[0,\infty)$
satisfying the Dini condition
\begin{equation}\label{eq:dini-condition}
 \int_0^1\frac{\omega(t)}{t}\,\rd t<\infty
\end{equation}
and
\begin{equation}\label{et2h0}
 h(t)\le Ct^4\omega(t),\qquad 0<t\le t_0.
\end{equation}
Then there exist $g_1,\ldots,g_n\in H^\infty(\mathbb B)$ such that
$$f=\sum_{\ell=1}^n f_\ell\ast g_\ell.$$

\item\label{gauge:slice}
Suppose that all  assumptions of \textup{(i)} hold with \eqref{et2h0} replaced by
$$h(t)\le Ct^2\omega(t),\qquad 0<t\le t_0.$$
If $f_1,\ldots,f_n$ preserve a common slice, then the conclusion of
\textup{(i)} holds.
\end{enumerate}

Moreover, in \textup{(i)} and \textup{(ii)}, if $f,f_1,\ldots,f_n\in H_+^\infty(\mathbb B)$, then the
coefficients $g_1,\ldots,g_n$ can also be chosen in $H_+^\infty(\mathbb B)$.
\end{thm}

The proof of Theorem~\ref{thm2}\textup{(i)} combines estimates on
quaternionic conjugacy spheres with scalar and matrix ideal methods.
Holloway and Trent~\cite[Sections~4.1 and~4.4]{HT} obtain ideal estimates
for matrix equations under determinant bounds, including matrix-valued
right-hand sides and concatenated generators. Here the sphere estimates
give the required control by the maximal minors (Lemma~\ref{ltrowestimate}).
Treil's Dini ideal theorem~\cite{Tr07} and an adjugate formula then yield
a bounded holomorphic solution, which is symmetrized to recover slice
regular coefficients. The power conclusions in
\cite[Section~4.3]{HT} use entrywise powers. In the present setting,
slice powers correspond to ordinary matrix powers:
\[
 \chi\left(f^{\ast N}\right)=\chi(f)^N.
\]

Taking $h(t)=Ct^p$ in Theorem~\ref{thm2} gives the following result.

\begin{cor}\label{thm:power-ideal}
Let $n\ge1$ be an integer, let $p\ge1$, and suppose that
$f,f_1,\ldots,f_n\in H^\infty(\mathbb B)$ satisfy
\[
 |f(q)|\le C\left(\sum_{\ell=1}^n\left|f_\ell(q)\right|\right)^p,
 \qquad q\in\mathbb B,
\]
for some constant $C>0$.
\begin{enumerate}
\renewcommand{\labelenumi}{\textup{(\roman{enumi})}}
\item
If $p>4$, then $f\in\mathcal I\left(f_1,\ldots,f_n\right)$.
\item
If $p>2$ and the generators preserve a common slice,
then
\[
 f\in\mathcal I\left(f_1,\ldots,f_n\right).
\]
\end{enumerate}
In \textup{(i)} and \textup{(ii)}, if $f,f_1,\ldots,f_n$ are intrinsic,  then the coefficients can also be chosen intrinsic.

For each $n\ge2$, there also exist
$f,f_1,\ldots,f_n\in H_+^\infty(\mathbb B)$ such that
\[
 \sum_{\ell=1}^n\left|f_\ell(q)\right|\le1,
 \qquad
 |f(q)|\le\left(\sum_{\ell=1}^n\left|f_\ell(q)\right|\right)^p
 \quad(q\in\mathbb B,\ 1\le p\le2),
\]
but $f\notin\mathcal I\left(f_1,\ldots,f_n\right)$.
\end{cor}

\begin{rem}
 For every $\varepsilon>0$,
Theorem~\ref{thm2}\textup{(i)} applies when, for
sufficiently small $t>0$,
\[
 h(t)\le C t^{4+\varepsilon},
 \qquad\text{or}\qquad
 h(t)\le\frac{C t^4}{[\log(e/t)]^{1+\varepsilon}},
\]
or
\[
 h(t)\le
 \frac{C t^4}
 {\log(e/t)\,[\log(e\log(e/t))]^{1+\varepsilon}}.
\]
For generators preserving a common slice, replace $t^{4+\varepsilon}$ by $t^{2+\varepsilon}$ and the numerators
$t^4$ by $t^2$.
 For related logarithmic conditions in the scalar setting,
see \cite[Remark~0.3]{Tr07}.
\end{rem}

\begin{rem}
 For $n\ge2$, the condition $p>2$ in Corollary~\ref{thm:power-ideal}\textup{(ii)} is sharp: failure at $p=2$ already occurs for intrinsic data.
 For general generators and $2<p\le4$, a remaining question is whether
\[
 |f|\le CS^p\quad\Longrightarrow\quad
 f\in\mathcal I(f_1,\ldots,f_n).
\]
One can also ask which weaker conditions on $\omega$ suffice in
Theorem~\ref{thm2}.

These bounds concern the pointwise quantities $|f|$ and $S$.
For nonintrinsic $f$, one need not have $|f^{\ast N}|=|f|^N$ when $N\ge2$.
\end{rem}

  Section~\ref{s2} collects the slice regular preliminaries and matrix
identities. Section~\ref{sec:counterexamples} proves Theorems \ref{t1}, \ref{thm2-1} and Corollary \ref{t2}. The sufficient conditions in Theorem \ref{thm2} and Corollary \ref{thm:power-ideal} are proved in Section~\ref{sec:proof-gauge}.

\section{Preliminaries}\label{s2}
 We recall the slice regular preliminaries and introduce the complex matrix representations used below. Our conventions follow \cite{GSS}; see also \cite{CPSW,GSaS,Lin}.  A superscript $*$ on a complex matrix denotes its conjugate transpose, and $\ast$ between functions denotes the slice product.

\subsection{Slice regular functions}
Denote the real  algebra of quaternions by
 \[
\mathbb H=\mathbb R+\mathbb Ri+\mathbb Rj+\mathbb Rk,
\]
where $i^2=j^2=k^2=ijk=-1$. For
$q=x_0+x_1i+x_2j+x_3k\in \mathbb{H}$ with $x_0,x_1,x_2,x_3\in\mathbb R$, define its conjugate and norm by
\[
 \bar q=x_0-x_1i-x_2j-x_3k, \qquad
 |q|=\left(q\bar q\right)^{1/2}  =\left(x_0^2+x_1^2+x_2^2+x_3^2\right)^{1/2}.
\]
In particular, $|ab|=|a||b|$ and $q^{-1}=\bar q/|q|^2$ for $q\ne0$.
Denote the quaternionic unit ball by
\[
\mathbb B=\{q\in \mathbb{H}: |q|<1\}.
\]
For $I$ in the set of imaginary units $\mathbb S=\left\{q\in\mathbb H:q^2=-1\right\}$, let $\mathbb C_I=\mathbb R+I\mathbb R$ be the corresponding complex plane and put
\[\mathbb D_I=\mathbb B\cap\mathbb C_I,\]
which  is the unit disk in $\mathbb C_I$.  We identify $\mathbb C_I$
with $\mathbb C$ by $x+Iy\longleftrightarrow x+iy$; this also identifies $\mathbb D_I$ with $\mathbb D$.
We write $H^\infty\left(\mathbb D_I\right)$ for the bounded
$\mathbb C_I$-valued holomorphic functions on $\mathbb D_I$,
with the supremum norm over that disk.
The conjugacy sphere of $z=x+Iy$ is
\[
[z]=\{x+Ky:K\in\mathbb S\},
\]
which reduces to a point if $y=0$. Every point of $[z]$ has the form $s^{-1}zs$ with $|s|=1$.

\begin{definition}
A   function $f:\mathbb B\to\mathbb H$ is called \emph{left slice regular} if, for every $I\in\mathbb S$, its
restriction $f_I=f|_{\mathbb D_I}$ has continuous partial derivatives and
\[ \frac12\left(\frac{\partial}{\partial x}   +I\frac{\partial}{\partial y}\right)f_I(x+Iy)=0.\]
\end{definition}
Throughout this paper, slice regular means left slice regular, also called
slice hyperholomorphic.
On $\mathbb B$, this is equivalent to the existence of a power
series
\begin{equation}\label{eqslseries}
 f(q)=\sum_{\nu=0}^{\infty}q^\nu a_\nu,
 \qquad a_\nu\in\mathbb H,
\end{equation}
converging locally uniformly on $\mathbb B$; see \cite[Theorems~1.6 and~1.10]{GSS}. Coefficients are always written on the right. The identity principle  \cite[Theorem~1.13]{GSS}  implies that
two slice regular functions agreeing on one slice, or on a nondegenerate real interval, agree throughout $\mathbb B$.

We use the following splitting and extension lemmas to pass between slice regular functions and holomorphic functions on a fixed slice.
\begin{lem}[Splitting { {\cite[Lemma~1.3, p.~2]{GSS}}}]
\label{lsplitting}
Fix orthogonal imaginary units $I,J\in\mathbb S$, so $IJ=-JI$.
Every slice regular function $f:\mathbb B\to\mathbb H$ has a unique splitting
\[
 f_I(z)=F(z)+G(z)J, \qquad F,G:\mathbb D_I\to\mathbb C_I\ \text{holomorphic}.
\]
For $f\in H^\infty(\mathbb B)$, we have
\begin{equation}\label{eq:splitting-bounds}
 \left|f_I(z)\right|^2=|F(z)|^2+|G(z)|^2,
 \qquad \max\left\{\|F\|_\infty,\|G\|_\infty\right\}\le\|f\|_\infty.
\end{equation}
\end{lem}

\begin{lem}[Extension   {\cite[Lemma~1.21, p.~10]{GSS}}]
\label{lextension}
Fix orthogonal imaginary units $I,J\in\mathbb S$, so $IJ=-JI$.
Let $F,G:\mathbb D_I\to\mathbb C_I$ be holomorphic, and put
$u_I=F+GJ$. Then $u_I$ has a unique slice regular extension
$\operatorname{ext}\left(u_I\right)$ to $\mathbb B$, given by
\begin{equation}\label{erep}
\operatorname{ext}\left(u_I\right)(x+Ky) =\frac{1-KI}{2}\,u_I(x+Iy)  +\frac{1+KI}{2}\,u_I(x-Iy),
\end{equation}
where $K\in\mathbb S$, $x,y\in\mathbb R$, and $x^2+y^2<1$.
If $u_I$ is bounded on $\mathbb D_I$, then
\[
 \left\|\operatorname{ext}\left(u_I\right)\right\|_\infty
 \le 2\sup_{z\in\mathbb D_I}\left|u_I(z)\right|.
\]
\end{lem}

For a $\mathbb{C}_I$-valued function $F$ on $\mathbb D_I$, put
\begin{equation}\label{eq:hat-involution}
 \widehat F(z)=\overline{F\left(\bar z\right)}, \quad z\in\mathbb D_I,
\end{equation}
where the bars denote conjugation in $\mathbb C_I$.
This is a conjugate-linear, multiplicative, isometric involution of $H^\infty\left(\mathbb D_I\right)$, i.e.
\[
\widehat{cF}=\bar c\,\widehat F, \quad c\in\mathbb C_I,
\]
and
\[
\widehat{FG}=\widehat F\,\widehat G,\ \  \|\widehat F \|_\infty=\|F\|_\infty, \ \ \widehat{\widehat F}=F.
\]

A function $F\in H^\infty\left(\mathbb D_I\right)$ is called \emph{real-symmetric} if $\widehat F=F$ on $\mathbb D_I$.
For matrix-valued functions, the hat is applied entrywise; it
does not include transposition.

Let $f$ and $g$ be slice regular functions on $\mathbb{B}$ with power series expansions
\begin{equation*}
 f(q)=\sum_{\nu=0}^{\infty}q^\nu a_\nu,\quad  g(q)=\sum_{\nu=0}^{\infty}q^\nu b_\nu.
\end{equation*}
Their slice product (or  $\ast$-product) \cite[Definition~1.36]{GSS}  is defined by
\begin{equation}\label{eqsliconvol}
 (f\ast g)(q)=\sum_{\nu=0}^{\infty}q^\nu
       \sum_{\mu=0}^{\nu}a_\mu b_{\nu-\mu}.
\end{equation}
For splittings $f_I=F+GJ$ and $g_I=U+VJ$, this gives
\begin{equation}\label{eqspprod}
 (f\ast g)_I=\left(FU-G\widehat V\right)+\left(FV+G\widehat U\right)J,
\end{equation}
which   agrees with the product used in \cite[Definition~1.37, p.~22]{GSS}.

 The product formula in \cite[Theorem~3.4, p.~42]{GSS}  gives the Banach algebra structure and norm estimate below.

\begin{lem}
\label{lbaal}
Let $f,g$ be slice regular on $\mathbb B$. Then, for every
$q\in\mathbb B$,
\begin{equation}\label{eqpstar}
 (f\ast g)(q)=
 \begin{cases}
 f(q)g\left(f(q)^{-1}qf(q)\right),&f(q)\ne0,\\
 0,&f(q)=0.
 \end{cases}
\end{equation}
Consequently, $\left( H^\infty(\mathbb B),+,\ast,\|\cdot\|_\infty\right)$ is a
unital real Banach algebra, and
\begin{equation}\label{eqfastg}
 \|f\ast g\|_\infty\le\|f\|_\infty\|g\|_\infty,
 \qquad f,g\in H^\infty(\mathbb B).
\end{equation}
\end{lem}

A slice regular function $f$ on $\mathbb B$ is \emph{intrinsic} if $f\left(\bar q\right)=\overline{f(q)}$ for all $q\in \mathbb B$. Equivalently, all the coefficients in \eqref{eqslseries}
are real, or $f\left(\mathbb D_I\right)\subseteq\mathbb C_I$ for every $I\in\mathbb S$; see \cite[Remarks~1.41 and~1.44]{GSS}.
 By \cite[Lemma~1.43]{GSS}, for intrinsic $f$ and arbitrary slice regular $g$,
\begin{equation}\label{eq:intrinsic-central}
f \ast g=g\ast f=fg,
\end{equation}
where  $fg$ is the pointwise product with $f$ on the left.
In particular, $f^{\ast N}(q)=f(q)^N$ for all $N\in \mathbb N$.

 For generators preserving a common slice, the ideal equation splits
into two scalar equations.

We denote the real-symmetric subalgebra by
\[
 H^\infty_{\mathbb R}(\mathbb D)
 =\left\{F\in H^\infty(\mathbb D):
          F\left(\bar z\right)=\overline{F(z)}
          \text{ for every }z\in\mathbb D\right\}.
\]

\begin{lem}\label{lintrinreduc}
Fix $I,J\in\mathbb S$ with $IJ=-JI$,  and identify $\mathbb D_I$ with $\mathbb D$.

(i) The restriction map
\[ H_+^\infty(\mathbb B) \longrightarrow H^\infty_{\mathbb R}(\mathbb D),
 \qquad f\longmapsto f_I=f|_{\mathbb D_I},\]
is an isometric isomorphism of real Banach algebras.

(ii) Let $f,f_1,\ldots,f_n\in H^\infty(\mathbb B)$, and suppose that
 \begin{equation}\label{one-slice}
f_\ell\left(\mathbb D_I\right)\subseteq\mathbb C_I, \qquad 1\le\ell\le n.
\end{equation}
Use $f_\ell$ also for its holomorphic restriction to
$\mathbb D_I$, and write
\[
 f_I=U+VJ,
 \qquad U,V\in H^\infty(\mathbb D).
\]
Then
\begin{equation}\label{ecidealred}
 f\in\mathcal I\left(f_1,\ldots,f_n\right)
 \quad\Longleftrightarrow\quad
 U,V\in\sum_{\ell=1}^n f_\ell H^\infty(\mathbb D).
\end{equation}

 (iii) If $f,f_1,\ldots,f_n\in H_+^\infty(\mathbb B)$ and
$f\in\mathcal I(f_1,\ldots,f_n)$, then there exist
$g_1,\ldots,g_n\in H_+^\infty(\mathbb B)$ such that
\[
 f=\sum_{\ell=1}^n f_\ell\ast g_\ell.
\]
\end{lem}

\begin{proof}(i)
For $F\in H^\infty(\mathbb D)$ with
$F(z)=\sum_{\nu\ge0}c_\nu z^\nu$, uniqueness of Taylor coefficients gives
\[
 \widehat F=F
 \quad\Longleftrightarrow\quad c_\nu\in\mathbb R\ \ (\nu\ge0)
 \quad\Longleftrightarrow\quad
 \operatorname{ext}(F)\text{ is intrinsic}.
\]
If $\widehat F=F$ and $F(x+Iy)=\alpha+I\beta$, with $\alpha,\beta\in\mathbb R$,
real symmetry and \eqref{erep} give
\[
 F(x-Iy)=\alpha-I\beta,\qquad
 \operatorname{ext}(F)(x+Ky)=\alpha+K\beta,
 \qquad K\in\mathbb S.
\]
Hence
\[
 |\operatorname{ext}(F)(x+Ky)|^2
 =\alpha^2+\beta^2=|F(x+Iy)|^2,
 \qquad
 \|\operatorname{ext}(F)\|_\infty=\|F\|_\infty.
\]
Restriction and extension are inverse real-linear maps by the
identity principle. For intrinsic $f,g$, \eqref{eq:intrinsic-central}
also gives $(f\ast g)_I=f_Ig_I$, proving the first assertion.

(ii) For the ideal reduction, write
$\left(g_\ell\right)_I=C_\ell+D_\ell J$ by Lemma~\ref{lsplitting}.
Since the generators satisfy (\ref{one-slice}), \eqref{eqspprod} gives
$$\Big(\sum_{\ell=1}^n f_\ell\ast g_\ell\Big)_I =\sum_{\ell=1}^n f_\ell C_\ell
  +\Big(\sum_{\ell=1}^n f_\ell D_\ell\Big)J.$$
Uniqueness of splitting therefore yields
\[
 f=\sum_{\ell=1}^n f_\ell\ast g_\ell
 \quad\Longrightarrow\quad
 U=\sum_{\ell=1}^n f_\ell C_\ell,\qquad
 V=\sum_{\ell=1}^n f_\ell D_\ell.
\]

Conversely, given these two scalar representations, set
\[
 g_\ell=\operatorname{ext}\left(C_\ell+D_\ell J\right)
       \in H^\infty(\mathbb B).
\]
Boundedness follows from Lemma~\ref{lextension}; the displayed
splitting identity gives
$\left(\sum_\ell f_\ell\ast g_\ell\right)_I\allowbreak=U+VJ=f_I$.
The identity principle proves \eqref{ecidealred}.

 (iii) For intrinsic data, part~(ii) gives $V=0$ and a scalar representation
of $U$. Since $\widehat U=U$ and $\widehat{f_\ell}=f_\ell$, reflecting this
representation and averaging the two expressions gives
\[
 U=\sum_{\ell=1}^n f_\ell C_\ell
  =\sum_{\ell=1}^n f_\ell\widehat C_\ell
  =\sum_{\ell=1}^n f_\ell\frac{C_\ell+\widehat C_\ell}{2}.
\]
 Each $(C_\ell+\widehat C_\ell)/2$ is bounded and real-symmetric.
By part~(i), its extension
$g_\ell=\operatorname{ext}((C_\ell+\widehat C_\ell)/2)$ is bounded and
intrinsic, and $f=\sum_\ell f_\ell\ast g_\ell$.
\end{proof}

\subsection{Complex matrix representations}

  For general generators, the two splitting components in
\eqref{eqspprod} are coupled. We write the ideal equation as a matrix
equation and express pointwise bounds on $\mathbb B$ in terms of row norms.

 Write $\mathbb C^{r\times s}$ for the complex $r\times s$ matrices,
$A_{\mu\nu}$ for the $(\mu,\nu)$ entry, and
$H^\infty(\mathbb D)^{r\times s}$ for matrices with entries in
$H^\infty(\mathbb D)$. Row and column vectors carry the Euclidean norm.
Denote
\[
 \|A\|_{\mathrm{op}}
 =\sup_{\substack{u\in\mathbb C^{1\times r}\\\|u\|=1}}\|uA\|
 =\sup_{\substack{v\in\mathbb C^s\\\|v\|=1}}\|Av\|,
 \qquad A\in\mathbb C^{r\times s},
\]
and use the associated supremum norm for matrix-valued functions.
We denote the $2\times2$ identity matrix by $I_2$ and use
$A^N$ for ordinary matrix powers.

Fix $I,J\in\mathbb S$ with $IJ=-JI$. For
$f\in H^\infty(\mathbb B)$ with splitting $f_I=F+GJ$,
its $2\times2$ complex matrix representation is given by
\begin{equation}\label{echidef}
 \chi(f)(z)=
 \begin{pmatrix}
 F(z)&G(z)\\-\widehat G(z)&\widehat F(z)
 \end{pmatrix}.
\end{equation}
Here $z\in\mathbb D_I$ is identified with a point of $\mathbb D$.
The pair $I,J$ is fixed throughout each argument. Related representations
appear in the quaternionic Wiener algebra setting in \cite{Alpay};
the map in \eqref{echidef} is denoted by $\mathcal R_I$ in
\cite{Lin}. The following lemma relates its row norms to values
on conjugacy spheres.

\begin{lem}[{\cite[Proposition~3.1]{Lin}}]
\label{lmatrep}
The map $\chi$ defined in \eqref{echidef} is an injective
real algebra homomorphism from $H^\infty(\mathbb B)$ into
the bounded holomorphic $2\times2$ matrix functions on
$\mathbb D$. In particular,
\begin{equation}\label{eq:chi-multiplicative}
 \chi(f\ast g)=\chi(f)\chi(g),
 \qquad f,g\in H^\infty(\mathbb B).
\end{equation}
For $f\in H^\infty(\mathbb B)$, $z\in\mathbb D_I$, and a unit row $u=(\alpha,\beta)\in\mathbb C_I^2$, set
$s=\alpha+\beta J$. Then
\begin{equation}\label{eq:row-sphere-identity}
 \|u\chi(f)(z)\|=\left|f\left(s^{-1}zs\right)\right|.
\end{equation}
Consequently,
\begin{equation}\label{eq:chi-isometry}
 \|\chi(f)(z)\|_{\mathrm{op}}=\sup_{q\in[z]}|f(q)|,
 \qquad
 \sup_{z\in\mathbb D}\|\chi(f)(z)\|_{\mathrm{op}}
 =\|f\|_\infty.
\end{equation}
\end{lem}

For $f,f_1,\ldots,f_n\in H^\infty(\mathbb B)$, set
\begin{equation}\label{eqblomat}
 H=\chi(f),\qquad
 M=\left[\,\chi\left(f_1\right)\ \cdots\ \chi\left(f_n\right)\,\right].
\end{equation}
The holomorphic matrices $H$ and $M$ have sizes $2\times2$
and $2\times2n$, respectively.

Fix $z\in\mathbb D_I$ and a unit row
$u=(\alpha,\beta)\in\mathbb C_I^2$, and set
\[
 s=\alpha+\beta J,\qquad q=s^{-1}zs.
\]
Then $|s|=1$ and $q\in[z]\subset\mathbb B$.
Equation \eqref{eq:row-sphere-identity} gives
(cf.\ \cite[(3.12)]{Lin})
\begin{equation}\label{eqbrnorm}
 \|uM(z)\|^2
 =\sum_{\ell=1}^n\left\|u\chi\left(f_\ell\right)(z)\right\|^2
  =\sum_{\ell=1}^n\left|f_\ell(q)\right|^2,\qquad
 \|uH(z)\|
 =\|u\chi(f)(z)\|=|f(q)|.
\end{equation}

For Hermitian matrices, $A\le B$ means that $B-A$ is
positive semidefinite. For each fixed $z$, \eqref{eqbrnorm} gives
\[
 \begin{aligned}
 H(z)H(z)^*\le M(z)M(z)^*
 &\quad\Longleftrightarrow\quad
 \|uH(z)\|^2\le\|uM(z)\|^2\quad\text{for every unit row }u\\
 &\quad\Longleftrightarrow\quad
 |f(q)|^2\le\sum_{\ell=1}^n\left|f_\ell(q)\right|^2
       \quad\text{for every }q\in[z].
 \end{aligned}
\]
Since $\mathbb B=\bigcup_{z\in\mathbb D_I}[z]$, we obtain
\begin{equation}\label{eqmorder}
 |f(q)|^2\le\sum_{\ell=1}^n\left|f_\ell(q)\right|^2
       \quad \forall q\in\mathbb B \qquad\Longleftrightarrow\quad
 H(z)H(z)^*\le M(z)M(z)^*
       \quad \forall z\in\mathbb D.
\end{equation}

More generally, suppose that
$h:[0,\infty)\to[0,\infty)$ is nondecreasing and
\[
 |f(q)|\le
 h\left(\sum_{\ell=1}^n\left|f_\ell(q)\right|\right),
 \qquad q\in\mathbb B.
\]
For the same $u,z,q$ as above, the Cauchy--Schwarz
inequality and \eqref{eqbrnorm} give
\[
 \sum_{\ell=1}^n\left|f_\ell(q)\right|
 \le\sqrt n
       \left(\sum_{\ell=1}^n\left|f_\ell(q)\right|^2\right)^{1/2}
 =\sqrt n\,\|uM(z)\|.
\]
By monotonicity of $h$, for every $z\in\mathbb D$ and unit row $u$ we obtain
\begin{equation}\label{eqgrabound}
 \|uH(z)\|
 =|f(q)| \le h\left(\sum_{\ell=1}^n\left|f_\ell(q)\right|\right)
 \le h\left(\sqrt n\,\|uM(z)\|\right).
\end{equation}

  By \eqref{eq:chi-multiplicative}, the ideal equation is $MX=H$,
where $X$ has the blocks $\chi(g_1),\ldots,\chi(g_n)$ stacked vertically.
 To obtain the slice regular coefficients, it is enough to solve
$MX=H$ with bounded holomorphic matrix blocks and then symmetrize
those blocks. This symmetry method appears in the matrix corona constructions of
\cite{CPSWmatrix,Lin}. We include the finite-block argument for a general
right-hand side.

\begin{lem}
\label{lmsymmetry}
For a holomorphic $2\times2$ matrix function $X$ on
$\mathbb D$, define
\[
 J_0=\begin{pmatrix}0&1\\-1&0\end{pmatrix},\qquad
 \mathfrak S(X)=J_0\widehat XJ_0^{-1},
\]
where $J_0^{-1}=-J_0$ and
$\widehat X(z)=\overline{X\left(\bar z\right)}$ is entrywise reflection,
without transposition. A bounded holomorphic $2\times2$
matrix $X$ has the form $X=\chi(g)$ for some
$g\in H^\infty(\mathbb B)$ if and only if $\mathfrak S(X)=X$.

 Let $n\ge1$ and $f,f_1,\ldots,f_n\in H^\infty(\mathbb B)$.
If bounded holomorphic $2\times2$ matrix functions
$X_1,\ldots,X_n$ on $\mathbb D$ satisfy
\[
 \sum_{\ell=1}^n\chi\left(f_{\ell}\right)X_{\ell}=\chi(f),
\]
 then $f\in\mathcal I\left(f_1,\ldots,f_n\right)$. One may take
\[
 \widetilde X_\ell=\frac{X_\ell+\mathfrak S\left(X_\ell\right)}2,
 \qquad
 g_\ell=\operatorname{ext}\left(\left(\widetilde X_\ell\right)_{11}
                         +\left(\widetilde X_\ell\right)_{12}J\right),
 \qquad 1\le\ell\le n.
\]
\end{lem}

\begin{proof}
Direct calculation shows that $\mathfrak S$ is a real-linear involution
and that
\[
 \mathfrak S(AB)=\mathfrak S(A)\mathfrak S(B).
\]
Its fixed points are precisely the matrices
\[
 X=\begin{pmatrix}F&G\\-\widehat G&\widehat F\end{pmatrix}.
\]
For bounded holomorphic $F,G$, this is $\chi(g)$ with
$g=\operatorname{ext}(F+GJ)\in H^\infty(\mathbb B)$.
Since $\chi(f)$ and each $\chi(f_\ell)$ are fixed by $\mathfrak S$,
applying the involution to the given equation yields
\[
 \sum_{\ell=1}^n\chi(f_\ell)\mathfrak S(X_\ell)=\chi(f).
\]
Averaging the two equations gives the same identity with
$X_\ell$ replaced by $\widetilde X_\ell$.
These matrices are bounded holomorphic fixed points, so
\eqref{eq:chi-multiplicative} gives the required right-ideal relation.
\end{proof}

The same symmetrization applies to Leech's factorization theorem
\cite{Leech}. It gives an operator criterion for the ideal equation,
which we will compare with pointwise domination.

For $r\ge1$, let $H^2(\mathbb D,\mathbb C^r)$ be the space of holomorphic
column vectors $\varphi(z)=\sum_{\nu\ge0}z^\nu\varphi_\nu$ with
\[
 \|\varphi\|_{H^2(\mathbb D,\mathbb C^r)}^2
 =\sum_{\nu\ge0}\|\varphi_\nu\|^2<\infty.
\]
For a bounded holomorphic $r\times s$ matrix function $A$, denote its multiplication operator by
\[
 \mathcal M_A:H^2\left(\mathbb D,\mathbb C^s\right)
       \longrightarrow H^2\left(\mathbb D,\mathbb C^r\right),
 \qquad \mathcal M_A\varphi=A\varphi.
\]
For operators, $^*$ denotes the Hilbert-space adjoint, and
inequalities are understood in the quadratic-form sense.
The following finite-family criterion follows from the classical Leech
theorem and Lemma~\ref{lmsymmetry}.

\begin{lem}[Leech factorization]\label{lqleech}
Let $n\ge1$ be an integer and $\gamma>0$. For
$f,f_1,\ldots,f_n\in H^\infty(\mathbb B)$, define $H$ and $M$
by \eqref{eqblomat}.
Then the following statements are equivalent.
\begin{enumerate}
\renewcommand{\labelenumi}{\textup{(\roman{enumi})}}
\item
The inequality
\begin{equation}\label{eq:leech-operator}
 \mathcal M_H\mathcal M_H^*
 \le\gamma^2\mathcal M_M\mathcal M_M^*
\end{equation}
holds on $H^2\left(\mathbb D,\mathbb C^2\right)$.
\item
There exist $g_1,\ldots,g_n\in H^\infty(\mathbb B)$ such that
\[
 f=\sum_{\ell=1}^n f_\ell\ast g_\ell,
 \qquad
 \sup_{z\in\mathbb D}
 \left\|\begin{pmatrix}
 \chi\left(g_1\right)(z)\\\vdots\\\chi\left(g_n\right)(z)
 \end{pmatrix}\right\|_{\mathrm{op}}\le\gamma.
\]
\end{enumerate}
In particular, $f\in\mathcal I\left(f_1,\ldots,f_n\right)$ if and only if
\eqref{eq:leech-operator} holds for some $\gamma>0$.
\end{lem}

\begin{proof}
Suppose \textup{(i)} holds. Leech's theorem~\cite{Leech} gives a bounded
holomorphic $2n\times2$ matrix $X$, with $2\times2$ blocks $X_\ell$,
such that $MX=H$ and $\|X\|_\infty\le\gamma$, where the norm is the
supremum of the pointwise operator norm. Put
\[
 J_n=\operatorname{diag}(J_0,\ldots,J_0),\qquad
 \widetilde X=\tfrac12\left(X+J_n\widehat XJ_0^{-1}\right).
\]
Its blocks are the symmetrized matrices in Lemma~\ref{lmsymmetry},
so $M\widetilde X=H$ and each block has the form $\chi(g_\ell)$.
Reflection and multiplication by the constant unitary matrices
$J_n,J_0^{-1}$ preserve the supremum norm. Hence
$\|\widetilde X\|_\infty\le\|X\|_\infty\le\gamma$, proving \textup{(ii)}.

Conversely, let $X$ be the matrix formed by stacking the blocks
$\chi(g_\ell)$ in \textup{(ii)}. Then $MX=H$ and
$\|\mathcal M_X\|=\|X\|_\infty\le\gamma$, so
\[
 \mathcal M_H\mathcal M_H^*
 =\mathcal M_M\mathcal M_X\mathcal M_X^*\mathcal M_M^*
 \le\gamma^2\mathcal M_M\mathcal M_M^*.
\]
\end{proof}

 To compare the operator inequality \eqref{eq:leech-operator} with
the pointwise condition \eqref{eqmorder}, we test it on reproducing kernels.

For $z,w\in\mathbb D$, set $k_z(w)=\left(1-w\bar z\right)^{-1}$.
For a bounded holomorphic $r\times s$ matrix $A$, the
reproducing-kernel identity gives
 \[
 \begin{aligned}
 \mathcal M_A^*(k_zv)&=k_zA(z)^*v,\\
 \|\mathcal M_A^*(k_zv)\|_{H^2(\mathbb D,\mathbb C^s)}^2
 &=\frac{v^*A(z)A(z)^*v}{1-|z|^2},\qquad v\in\mathbb C^r.
 \end{aligned}
\]
 Testing \eqref{eq:leech-operator} on $k_zv$, $v\in\mathbb C^2$,
and cancelling $\left(1-|z|^2\right)^{-1}$ gives
\[
 H(z)H(z)^*\le\gamma^2M(z)M(z)^*,\qquad z\in\mathbb D.
\]
 For the data in Theorem~\ref{t1}, the pointwise inequality holds with
$\gamma=1$, while \eqref{eq:leech-operator} fails for every $\gamma>0$
because $f\notin\mathcal I\left(f_1,f_2\right)$.

\smallskip
\noindent\emph{The single-generator case.}
Pointwise domination implies right-ideal membership for a single
generator. Indeed, suppose that
$f,f_1\in H^\infty(\mathbb B)$ and
$|f(q)|\le C\left|f_1(q)\right|$ on $\mathbb B$, where $C\ge0$.
If $f_1\equiv0$, then $f=0$. Otherwise, set
$A=\chi\left(f_1\right)$ and $H=\chi(f)$. The row-sphere identity gives
\[
 HH^*\le C^2AA^*.
\]
For real $x\in(-1,1)$, we have
$\det A(x)=\left|f_1(x)\right|^2$. Since $f_1\not\equiv0$, the identity
principle shows that $f_1$ cannot vanish on all of $(-1,1)$.
Thus $\det A\not\equiv0$. Away from its isolated zeros,
the holomorphic matrix $X=A^{-1}H$ satisfies
\[
 XX^*=A^{-1}HH^*\left(A^*\right)^{-1}\le C^2I_2,
 \qquad \|X\|_{\mathrm{op}}\le C.
\]
Thus $X$ extends to a bounded holomorphic matrix on
$\mathbb D$ by the removable singularity theorem, and
$AX=H$ holds everywhere. Lemma~\ref{lmsymmetry} yields
$f\in\mathcal I\left(f_1\right)$.

\section{Counterexamples}\label{sec:counterexamples}

In this section, we use two explicit Blaschke products to prove Theorem~\ref{t1}.
For the quadratic obstruction in Theorem \ref{thm2-1}, we localize Treil's scalar counterexample and impose real symmetry.

\subsection{Failure for all positive slice powers}\label{sec:proof-all-powers}

A prime denotes differentiation with respect to the complex
variable; for matrix functions, derivatives are taken entrywise.
 By the Schwarz--Pick lemma \cite{G}, if
$V\in H^\infty(\mathbb D)$, then
\begin{equation}\label{eq:bounded-derivative}
 \left(1-|z|^2\right)\left|V'(z)\right|\le\|V\|_\infty,
 \qquad z\in\mathbb D.
\end{equation}
 The next lemma gives the
Blaschke product estimates used in the counterexample;
see \cite{G} for the standard product formulas.

\begin{lem}
\label{lem:blaschke-sequences}
For $k\ge1$, let
\[
 r_k=1-2^{-k},\qquad
 s_k=r_k+\frac{\left(1-r_k\right)^2}{4},\qquad z_k=ir_k,
\]
and define the Blaschke products
\begin{equation}\label{eq:two-blaschke-products}
 B(z)=\prod_{k=1}^{\infty}\frac{r_k+iz}{1+ir_kz},
 \qquad
 C(z)=\prod_{k=1}^{\infty}\frac{s_k+iz}{1+is_kz}.
\end{equation}
Then there exist constants
$\delta_0,\delta_1>0$ such that, for every $k\ge1$,
\begin{align}
 \left(1-\left|z_k\right|^2\right)\left|B'\left(z_k\right)\right|&\ge\delta_0,
 \label{eq:blaschke-derivative-lower}\\
 \left|\widehat B\left(z_k\right)\right|&\ge\delta_1,
 \label{eq:reflected-blaschke-lower}\\
 0<\left|C\left(z_k\right)\right|&\le\frac{1-r_k}{4}\longrightarrow0.
 \label{eq:nearby-blaschke-small}
\end{align}
\end{lem}

\begin{proof}
The sums $\sum_k\left(1-r_k\right)$ and $\sum_k\left(1-s_k\right)$ are finite,
and $r_k<s_k<r_{k+1}$. Hence both products converge, their
zeros are simple, and their zero sets are disjoint.
For $k\ne\ell$, comparison of $1-r_k=2^{-k}$ and
$1-r_\ell=2^{-\ell}$ gives
\[
 \frac{\left|r_k-r_\ell\right|}{1-r_kr_\ell}
 \ge\frac{1-2^{-|k-\ell|}}{1+2^{-|k-\ell|}}.
\]
The derivative formula at a simple Blaschke zero therefore gives
\[
 \left(1-r_k^2\right)\left|B'\left(ir_k\right)\right|
 =\prod_{\ell\ne k}\frac{\left|r_k-r_\ell\right|}{1-r_kr_\ell}
 \ge\prod_{m=1}^{\infty}
       \left(\frac{1-2^{-m}}{1+2^{-m}}\right)^2
 =:\delta_0>0.
\]
The product is positive because
\[
 \sum_{m\ge1}\left[1-
       \left(\frac{1-2^{-m}}{1+2^{-m}}\right)^2\right]
 =\sum_{m\ge1}\frac{4\cdot2^{-m}}{\left(1+2^{-m}\right)^2}<\infty.
\]
Moreover,
\[
 \left|\widehat B\left(ir_k\right)\right|
 =\prod_{\ell=1}^{\infty}
       \frac{r_\ell+r_k}{1+r_\ell r_k}
 \ge\prod_{\ell=1}^{\infty}r_\ell
 =:\delta_1>0.
\]
All factors of $C$ have modulus at most one in the disk,
so the factor with index $k$ gives
\[
 \left|C\left(ir_k\right)\right|\le\frac{s_k-r_k}{1-s_kr_k}
 \le\frac{\left(1-r_k\right)^2}{4\left(1-r_k\right)}
 =\frac{1-r_k}{4}.
\]
It is nonzero because $ir_k$ is not a zero of $C$.
\end{proof}

\begin{proof}[Proof of Theorem~\ref{t1}]
Fix the standard orthogonal imaginary units $i,j$, and let
$B,C$ and $z_k=ir_k$ be as in
Lemma~\ref{lem:blaschke-sequences}. Set
\begin{equation}\label{eq:counterexample-quaternion-data}
 f_1=\operatorname{ext}\left(B^2\right),\qquad
 f_2=\operatorname{ext}\left(C^2\right),\qquad
 f=\frac12\operatorname{ext}\left(B^2+BCj\right).
\end{equation}
By Lemma~\ref{lextension}, these functions belong to
$H^\infty(\mathbb B)$. Their matrix representations are
\begin{equation}\label{eq:counterexample-matrices}
 \chi\left(f_1\right)=\begin{pmatrix}B^2&0\\0&\widehat B^{\,2}\end{pmatrix},
 \qquad
 \chi\left(f_2\right)=\begin{pmatrix}C^2&0\\0&\widehat C^{\,2}\end{pmatrix},\qquad
 H=\chi(f)=\frac12
 \begin{pmatrix}B^2&BC\\-\widehat B\,\widehat C&\widehat B^{\,2}\end{pmatrix}.
\end{equation}

Now let us prove
$HH^*\le\sum_{\ell=1}^2\chi(f_\ell)\chi(f_\ell)^*$.
Fix $z\in\mathbb D$ and put
\[
 d_1=|B(z)|^4+|C(z)|^4,\qquad
 d_2=\left|\widehat B(z)\right|^4+\left|\widehat C(z)\right|^4.
\]
Suppressing the common argument $z$, the inequality
$2|BC|^2\le|B|^4+|C|^4$ and its reflected version give
\[
 |B|^4+|BC|^2\le\tfrac32d_1,\qquad
 \left|\widehat B\right|^4+\left|\widehat B\widehat C\right|^2\le\tfrac32d_2.
\]
For every row $u=(\alpha,\beta)\in\mathbb C^2$, we therefore have
\begin{align*}
 \|uH\|^2
 &=\frac14\left(
     \left|\alpha B^2-\beta\widehat B\,\widehat C\right|^2
       +\left|\alpha BC+\beta\widehat B^{\,2}\right|^2\right)\\
 &\le\frac12\left(
      |\alpha|^2\left(|B|^4+|BC|^2\right)
       +|\beta|^2\left(\left|\widehat B\, \widehat C\right|^2+\left|\widehat B\right|^4\right)\right)
 \le\frac34\left(|\alpha|^2d_1+|\beta|^2d_2\right).
\end{align*}
Consequently,
\[
 HH^*\le\frac34
    \begin{pmatrix}d_1&0\\0&d_2\end{pmatrix}
 =\frac34\left(\chi\left(f_1\right)\chi\left(f_1\right)^*+\chi\left(f_2\right)\chi\left(f_2\right)^*\right).
\]
Hence, by \eqref{eqmorder},
\begin{equation}\label{eq:counterexample-domination}
 |f(q)|^2\le\frac34\left(\left|f_1(q)\right|^2+\left|f_2(q)\right|^2\right)
 \le \left|f_1(q)\right|^2+\left|f_2(q)\right|^2,
 \qquad q\in\mathbb B.
\end{equation}

 It remains to exclude every positive slice power from the ideal.
 Fix an integer $N\ge1$ and suppose that
$f^{\ast N}=f_1\ast g_1+f_2\ast g_2$ for some
$g_1,g_2\in H^\infty(\mathbb B)$.  Write
$\left(g_\ell\right)_i=U_\ell+V_\ell j$, where
$U_\ell,V_\ell\in H^\infty(\mathbb D)$ for $\ell=1,2$.
Put
\[
 K_N(z)=\left(H(z)^N\right)_{12}.
\]
Since $\chi\left(f^{\ast N}\right)=H^N$, comparison of the upper right
entries gives
\begin{equation}\label{eq:counterexample-coefficient-equation}
 K_N=B^2V_1+C^2V_2.
\end{equation}
Evaluate this identity and its derivative at the zeros $z_k$ of $B$.

 At $z_k$, we have $B(z_k)=0$, $C(z_k)\ne0$, and
\[
 H\left(z_k\right)=\frac12
 \begin{pmatrix}0&0\\
 -\widehat B\left(z_k\right)\widehat C\left(z_k\right)&\widehat B\left(z_k\right)^2\end{pmatrix}.
\]
For $m\ge1$, the first row of $H\left(z_k\right)^m$ is zero. Hence
\begin{equation}\label{V-zero}
 K_N\left(z_k\right)=0=C(z_k)^2V_2\left(z_k\right).\end{equation}

By (\ref{eq:nearby-blaschke-small}), we have $C(z_k)\neq0$ and then   $V_2\left(z_k\right)=0.$
Differentiating \eqref{eq:counterexample-coefficient-equation}
and using $B\left(z_k\right)=V_2\left(z_k\right)=0$ gives
\begin{equation}\label{eq:counterexample-derivative-equation}
 K_N'\left(z_k\right)
 =\left(2BB'V_1+B^2V_1'+2CC'V_2+C^2V_2'\right)\left(z_k\right) =C\left(z_k\right)^2V_2'\left(z_k\right).
\end{equation}
  To estimate $V_2'(z_k)$, we compute $K_N'(z_k)$ from the matrix
product rule
\[
 (H^N)'=\sum_{\ell=0}^{N-1}H^\ell H'H^{N-1-\ell}.
\]
At $z_k$, every term with $\ell\ge1$ has zero first row, so
$K_N'(z_k)=(H'(z_k)H(z_k)^{N-1})_{12}$.
The first row of $H'(z_k)$ and the second column of $H(z_k)^{N-1}$ are
\[
 \left(0,\frac{B'(z_k)C(z_k)}2\right),
 \qquad
 \begin{pmatrix}0\\(\widehat B(z_k)^2/2)^{N-1}\end{pmatrix},
\]
respectively, with $H^0=I_2$ when $N=1$. Their product gives
\begin{equation}\label{eq:counterexample-power-derivative}
 K_N'(z_k)=\frac{B'\left(z_k\right)C\left(z_k\right)\widehat B\left(z_k\right)^{2N-2}}{2^N}.
\end{equation}
By Lemma~\ref{lem:blaschke-sequences}, we have $K_N'(z_k)\neq0$. Together with \eqref{V-zero}, this shows that each $z_k$ is a simple zero of $K_N$. Combining \eqref{eq:counterexample-derivative-equation}
and \eqref{eq:counterexample-power-derivative} gives
\[
 V_2'\left(z_k\right)=
 \frac{B'\left(z_k\right)\widehat B\left(z_k\right)^{2N-2}}{2^NC\left(z_k\right)}.
\]
Lemma~\ref{lem:blaschke-sequences} and
\eqref{eq:bounded-derivative} give
\[
 \left\|V_2\right\|_\infty
 \ge\left(1-\left|z_k\right|^2\right)\left|V_2'\left(z_k\right)\right|
 \ge\frac{\delta_0\delta_1^{2N-2}}{2^N\left|C\left(z_k\right)\right|}
 \longrightarrow\infty,
\]
 contrary to $V_2\in H^\infty(\mathbb D)$.
The functions $f,f_1,f_2$ were chosen independently of $N$, so the
same triple has the asserted property for every positive integer $N$.
 \end{proof}

\subsection{A localized scalar counterexample}

 The following adaptation of \cite[Theorem~1.2 and Section~3]{Tr} gives a scalar counterexample whose three functions are uniformly bounded away from zero outside a prescribed boundary neighborhood.
The reflection argument in the next subsection uses these lower bounds.

\begin{lem}
\label{lloctreil}
Fix $0<\varepsilon<1/2$, and put
$U=\{z\in\mathbb D:|z-i|<\varepsilon\}$.
There exist $a,b,v\in H^\infty(\mathbb D)$, each of norm
at most one, and a constant $c>0$ such that
\begin{align}
 |v|&\le |a|+|b|\quad\text{on }\mathbb D,
 \label{eqlocdom}\\
 v^2&\notin aH^\infty(\mathbb D)+bH^\infty(\mathbb D),
 \label{eq:localized-nonmembership}\\
 \min\{|a(z)|,|b(z)|,|v(z)|\}&\ge c
       \quad(z\in\mathbb D\setminus U).
 \label{eq:localized-lower-bound}
\end{align}
\end{lem}

\begin{proof}
Let $\mathbb C_+=\{w\in\mathbb C:\operatorname{Im}w>0\}$
be the upper half-plane, and let $H^\infty\left(\mathbb C_+\right)$
denote its bounded holomorphic functions, with the supremum norm.
We construct finite interpolation problems with growing lower bounds
for solution norms, then combine widely separated translates into
infinite products. Keeping the factors close to one away from their
respective half-disks will give the lower bounds outside $U$.

\smallskip
\noindent\emph{Step 1. Finite interpolation and correction estimates.}
By \cite[Theorem~1.2]{Tr}, there is an
absolute constant $C_*>0$ such that, for each sufficiently small
$\delta>0$, there are finite Blaschke products $A_\delta,B_\delta$
on $\mathbb D$, with $B_\delta$ having only simple zeros, such that
\[
 \left|A_\delta\right|+\left|B_\delta\right|\ge\delta,\qquad
 A_\delta G_1+B_\delta G_2=1
 \ \Longrightarrow\
 \left\|G_1\right\|_\infty
 \ge C_*\delta^{-2}\log\log(1/\delta),
\]
where $G_1,G_2\in H^\infty(\mathbb D)$.
For each $k\ge1$, choose $0<\delta_k<e^{-1}$ sufficiently
small that
\[
 C_*\log\log\left(1/\delta_k\right)\ge3k.
\]
Transfer the two products to $\mathbb C_+$ by a Cayley map
and multiply them by unimodular constants so that their limits
at infinity are one. Denote the resulting products by
$A_k,B_k$, and the finite zero set of $B_k$ by $\sigma_k$.
These operations preserve both the corona lower bound and the
estimate for the first coefficient. Thus
\[
 \left\|A_k\right\|_\infty,\left\|B_k\right\|_\infty\le1,
 \qquad \left|A_k\right|+\left|B_k\right|\ge\delta_k.
\]
In particular, $A_k$ is nonzero on $\sigma_k$.

If $G\in H^\infty\left(\mathbb C_+\right)$ satisfies
$G(\lambda)=\delta_k^2/A_k(\lambda)$ on $\sigma_k$,  finite Blaschke division (see \cite[Section~2]{Tr}) gives
\[
 Q=\frac{\delta_k^2-A_kG}{B_k}\in H^\infty\left(\mathbb C_+\right),
 \qquad A_k\frac{G}{\delta_k^2}+B_k\frac{Q}{\delta_k^2}=1.
\]
Treil's estimate therefore gives
\[
 \|G\|_\infty
 \ge C_*\log\log\left(1/\delta_k\right)\ge3k.
\]
 We need this lower bound to survive small perturbations of the
prescribed values on $\sigma_k$. Since $\sigma_k$ is finite, the
interpolation functions below give an explicit bound for the correction.

For $\lambda\in\sigma_k$, choose
$\ell_{k,\lambda}\in H^\infty\left(\mathbb C_+\right)$ with
$\ell_{k,\lambda}(\lambda)=1$ and
$\ell_{k,\lambda}(\mu)=0$ for $\mu\in\sigma_k\setminus\{\lambda\}$.
Lagrange polynomials on the disk, composed with a Cayley map,
provide such functions. Put
\[
 L_k=\sum_{\lambda\in\sigma_k}\left\|\ell_{k,\lambda}\right\|_\infty.
\]
Given $e_\lambda\in\mathbb C$ for $\lambda\in\sigma_k$,
\[
 E=\sum_{\lambda\in\sigma_k}e_\lambda\ell_{k,\lambda}
 \quad\Longrightarrow\quad
 E(\lambda)=e_\lambda,\qquad
 \|E\|_\infty\le L_k\max_{\lambda\in\sigma_k}\left|e_\lambda\right|.
\]
\smallskip
\noindent\emph{Step 2. Local factors and stable interpolation.}
We replace the constant $\delta_k$ by a zero-free holomorphic function
that is close to $\delta_k$ on $\sigma_k$ and tends to one at infinity.

For $R>0$, define the zero-free holomorphic function
\begin{equation}\label{eqwkrw}
 W_{k,R}(w)=\exp\left(
    \frac{\log\delta_k}{\pi i}
    \int_{-R}^{R}\frac{\rd t}{t-w}\right),
 \qquad w\in\mathbb C_+.
\end{equation}
For $w=x+iy$, $y>0$, the real part of the exponent is
\[
 \frac{\log\delta_k}{\pi}
       \int_{-R}^{R}\frac{y\,\rd t}{(t-x)^2+y^2}.
\]
The integral divided by $\pi$ lies between zero and one;
hence $\delta_k\le\left|W_{k,R}\right|\le1$.
Moreover, for $|w|<R$ in the upper half-plane,
\[
 \int_{-R}^{R}\frac{\rd t}{t-w}
 =\pi i+\log(1-w/R)-\log(1+w/R),
\]
 where the logarithms on the right are the principal branches.
For any fixed $\rho>0$ and $R>\rho$, this expression extends
continuously to $\{|w|\le\rho,\ \operatorname{Im}w\ge0\}$.
Its two logarithmic terms tend uniformly to zero as $R\to\infty$.
Thus $W_{k,R}\to\delta_k$ uniformly on this closed half-disk,
with the values on the real axis understood as limits from above.  For fixed $k,R$ and $|w|>R$,
\[
 \left|\frac{\log\delta_k}{\pi i}
       \int_{-R}^{R}\frac{\rd t}{t-w}\right|
 \le\frac{2R\left|\log\delta_k\right|}{\pi(|w|-R)}
 \longrightarrow0\qquad(|w|\to\infty).
\]
 The estimate holds for $w\in\overline{\mathbb C_+}$ with $|w|>R$, where boundary values are taken outside $[-R,R]$. Thus $W_{k,R}(w)\to1$ uniformly as $|w|\to\infty$ in the closed upper half-plane.

For each fixed $k$, choose $r>0$ so large that $\sigma_k\subset\{|w|<r\}$ and
$\left|A_k(w)\right|+\left|B_k(w)\right|\ge1/2$ for $|w|\ge r$.
This is possible because both products tend to one at infinity.
By increasing $R>r$, we can arrange simultaneously that
$\left|W_{k,R}\right|\le2\delta_k$ on $|w|\le r$ and
\[
 \max_{\lambda\in\sigma_k}
 \left|\frac{W_{k,R}(\lambda)^2-\delta_k^2}
              {A_k(\lambda)}\right|
 L_k\le k.
\]
Fix these radii, which may depend on $k$, and put $W_k=W_{k,R}$.
The inner and outer estimates give

\[
 \left|W_k(w)\right|\le
 \begin{cases}
 2\delta_k\le2\left(\left|A_k(w)\right|+\left|B_k(w)\right|\right),&|w|\le r,\\
 1\le2\left(\left|A_k(w)\right|+\left|B_k(w)\right|\right),&|w|>r.
 \end{cases}
\]
Thus
\begin{equation}\label{eq:local-block-domination}
 \left|W_k(w)\right|\le2\left(\left|A_k(w)\right|+\left|B_k(w)\right|\right),
 \qquad w\in\mathbb C_+.
\end{equation}
 We next use the interpolation condition imposed on $W_k$ to retain
the lower bound for solution norms. Set
\[
 d_{k,\lambda}=\frac{W_k(\lambda)^2}{A_k(\lambda)},
 \qquad \lambda\in\sigma_k,
\]
and choose $\epsilon_k>0$ with $\epsilon_kL_k\le k$.
If $G\in H^\infty\left(\mathbb C_+\right)$ satisfies
\[
 \max_{\lambda\in\sigma_k}\left|G(\lambda)-d_{k,\lambda}\right|<\epsilon_k,
\]
define
\[
 G_0=G+\sum_{\lambda\in\sigma_k}
       \left(\frac{\delta_k^2}{A_k(\lambda)}-G(\lambda)\right)
          \ell_{k,\lambda}.
\]
Then $G_0(\lambda)=\delta_k^2/A_k(\lambda)$ on $\sigma_k$.
The choice of $R$ and the definition of $L_k$ give
 \begin{equation}\label{eq:stable-interpolation-bound}
 \|G_0-G\|_\infty\le k+L_k\epsilon_k\le2k,
 \qquad \|G\|_\infty\ge\|G_0\|_\infty-2k\ge k.
\end{equation}
 The other factors in the infinite products will introduce multiplicative
errors in $d_{k,\lambda}$. To keep the resulting absolute errors below
$\epsilon_k$, put
\[
 \theta_k=\min\left\{\frac14,
 \frac{\epsilon_k}
 {2\left(1+\max_{\lambda\in\sigma_k}\left|d_{k,\lambda}\right|\right)}\right\}.
\]

\smallskip
\noindent\emph{Step 3. Separated translates and infinite products.}
We translate the triples $(A_k,B_k,W_k)$ far enough apart that
the other factors change each prescribed value by less than $\epsilon_k$.
To obtain the lower bounds outside $U$, we also keep their half-disks
away from the region corresponding to $\mathbb D\setminus U$.
Let $\Phi:\mathbb C_+\to\mathbb D$ be the conformal bijection
$\Phi(w)=i(w-i)/(w+i)$, which sends infinity to $i$.  Since
$|\Phi(w)-i|=2/|w+i|$, the inverse image of
$\mathbb D\setminus U$ is contained in the compact set
\[
 K=\left\{w\in\overline{\mathbb C_+}:|w+i|\le2/\varepsilon\right\}.
\]
Choose $R_k>0$ so large that
the zeros of $A_k,B_k$ and the boundary interval associated
with $W_k$ lie in $|w|<R_k$, and
\begin{equation}\label{eq:local-block-tail}
 \max\left\{\left|A_k(w)-1\right|,\left|B_k(w)-1\right|,\left|W_k(w)-1\right|\right\}
 \le2^{-k-3}\qquad\left(|w|\ge R_k, w\in\overline{\mathbb C_+}\right).
\end{equation}
Values on the real axis here are taken outside that interval,
where the factors extend continuously.

We choose positive real translations $\tau_k$ recursively.
First require the half-disks
\[
 \Delta_k=\left\{w\in\overline{\mathbb C_+}:\left|w-\tau_k\right|\le R_k\right\}
\]
to be pairwise disjoint, disjoint from $K$, and to satisfy
$\tau_k-R_k\ge k$ for every $k\ge1$.  In particular,
 $\tau_k-R_k\to\infty$. We also control the earlier factors at the new
interpolation nodes and the new factor at all earlier nodes.
The translated zero set is
$\tau_k+\sigma_k=\{\tau_k+\lambda:\lambda\in\sigma_k\}$.
An empty product is understood to equal one. On this new set,
require
\begin{equation}\label{eqwwaw}
 \left|\prod_{\ell<k}
   \frac{W_{\ell}\left(w-\tau_{\ell}\right)^2}{A_{\ell}\left(w-\tau_{\ell}\right)}-1\right|
 \le\frac{\theta_k}{4}.
\end{equation}
On every earlier set $\tau_\ell+\sigma_\ell$, $\ell<k$,
require
\begin{equation}\label{eqwkwtk}
 \left|\frac{W_k\left(w-\tau_k\right)^2}{A_k\left(w-\tau_k\right)}-1\right|
 \le\theta_\ell 2^{-k-3}.
\end{equation}
 For fixed $k$, all these requirements hold when $\tau_k$ is sufficiently large. On the new finite set, every earlier argument $w-\tau_\ell$ tends to infinity; on each earlier finite set, $w-\tau_k$ tends to infinity in the opposite real direction. The corresponding factors tend uniformly to one, so the denominators are nonzero. Increasing $\tau_k$ also makes $\Delta_k$ disjoint from $K$ and the preceding half-disks.

 With these translations fixed, define
\[
 A(w)=\prod_{k\ge1}A_k\left(w-\tau_k\right),\qquad
 B_0(w)=\prod_{k\ge1}B_k\left(w-\tau_k\right),\qquad
 W(w)=\prod_{k\ge1}W_k\left(w-\tau_k\right).
\]
For every compact $E\subset\mathbb C_+$, the condition
$\tau_k-R_k\to\infty$ implies $E\cap\Delta_k=\varnothing$
for all $k\ge k_E$, with $k_E$ sufficiently large. Thus
\[
 \sum_{k\ge k_E}\sup_{w\in E}
 \max\left\{\left|A_k\left(w-\tau_k\right)-1\right|,\left|B_k\left(w-\tau_k\right)-1\right|,
                         \left|W_k\left(w-\tau_k\right)-1\right|\right\}
 \le\sum_{k\ge k_E}2^{-k-3}<\infty.
\]
The products converge locally uniformly. On each compact $E$,
their tails starting at $k_E$ have modulus at least
$\prod_{k\ge k_E}\left(1-2^{-k-3}\right)>0$.
Thus $A,B_0,W\in H^\infty\left(\mathbb C_+\right)$ have norms at most one,
and their zeros are precisely those of their factors. In particular,
\[
 W(w)\ne0\quad\left(w\in\mathbb C_+\right),\qquad
 A\left(\tau_k+\lambda\right)\ne0\quad\left(\lambda\in\sigma_k\right),
\]
because $A_k(\lambda)\ne0$ and the other factors obey
\eqref{eq:local-block-tail} on $\tau_k+\sigma_k$.
\smallskip
\noindent\emph{Step 4. Domination, localization, and nonmembership.}
We first pass \eqref{eq:local-block-domination} to the infinite products.
The tail estimates give the uniform lower bound
\[
 \prod_{k\ge1}\left(1-2^{-k-3}\right)\ge1-\sum_{k\ge1}2^{-k-3}
 =\frac78.
\]
  Together with \eqref{eq:local-block-domination} and
\eqref{eq:local-block-tail}, this gives, for $w\in\Delta_k\cap\mathbb C_+$,
\[
 |A(w)|+\left|B_0(w)\right|
 \ge\frac78\left(\left|A_k\left(w-\tau_k\right)\right|+\left|B_k\left(w-\tau_k\right)\right|\right) \ge\frac7{16}\left|W_k\left(w-\tau_k\right)\right|
 \ge\frac14|W(w)|.
\]
Outside all the half-disks, both $|A|$ and $\left|B_0\right|$ are
at least $7/8$, while $|W|\le1$, so the same final
inequality holds there. In particular,
\begin{equation}\label{eqwabc}
 |W|\le4\left(|A|+\left|B_0\right|\right)\quad\text{on }\mathbb C_+.
\end{equation}
 Since $K$ is disjoint from every $\Delta_k$, the tail estimate
\eqref{eq:local-block-tail} applies to every factor on $K\cap\mathbb C_+$.
Hence all three products have modulus at least $7/8$ there.
Under $\Phi$, this gives the required lower bounds outside $U$.

 To prove that $W^2$ is outside the ideal generated by $A,B_0$, we
estimate $W^2/A$ at the translated zeros of $B_k$.
For $w=\tau_k+\lambda$, $\lambda\in\sigma_k$, write
\[
 E_k(w)=\prod_{\ell<k}
           \frac{W_{\ell}\left(w-\tau_{\ell}\right)^2}{A_{\ell}\left(w-\tau_{\ell}\right)},\qquad
 T_k(w)=\prod_{\ell>k}
           \frac{W_{\ell}\left(w-\tau_{\ell}\right)^2}{A_{\ell}\left(w-\tau_{\ell}\right)}.
\]
 At these nodes, every denominator is nonzero. By
\eqref{eqwkwtk}, the deviations of the factors of $T_k$ from one
are summable, so this product converges. Passing to the limit in
the corresponding finite-product identity gives
\[
 \frac{W(\tau_k+\lambda)^2}{A(\tau_k+\lambda)}
 =d_{k,\lambda}\,E_k(\tau_k+\lambda)T_k(\tau_k+\lambda).
\]
We show that $|E_kT_k-1|<\theta_k$ at these nodes.
Equations \eqref{eqwwaw} and
\eqref{eqwkwtk} give
\[
 \left|E_k(w)-1\right|\le\frac{\theta_k}{4},\qquad
 \sum_{\ell>k}\left|
   \frac{W_\ell\left(w-\tau_\ell\right)^2}{A_\ell\left(w-\tau_\ell\right)}-1\right|
 \le\theta_k\sum_{\ell>k}2^{-\ell-3}\le\frac{\theta_k}{16}.
\]
Using $\left|\prod_\ell\left(1+u_\ell\right)-1\right|
\le\exp\left(\sum_\ell\left|u_\ell\right|\right)-1$ for summable $\left(u_\ell\right)$, we obtain
\[
 \left|T_k(w)-1\right|\le e^{\theta_k/16}-1\le\frac{\theta_k}{8}.
\]
Consequently,
\[
 \left|E_k(w)T_k(w)-1\right|
 \le\frac{\theta_k}{4}
       +\left(1+\frac{\theta_k}{4}\right)\frac{\theta_k}{8}
 <\theta_k.
\]
Hence, for $w=\tau_k+\lambda$,
\[
 \left|\frac{W(w)^2}{A(w)}-d_{k,\lambda}\right|
 =\left|d_{k,\lambda}\right|\,\left|E_k(w)T_k(w)-1\right|
 <\left|d_{k,\lambda}\right|\theta_k<\epsilon_k.
\]
 Now suppose that $W^2=AG_1+B_0G_2$, with
$G_1,G_2\in H^\infty(\mathbb C_+)$. Since $B_0$ vanishes on
$\tau_k+\sigma_k$, evaluation there gives
\[
 G_1\left(\tau_k+\lambda\right)
 =\frac{W\left(\tau_k+\lambda\right)^2}{A\left(\tau_k+\lambda\right)},
 \qquad \lambda\in\sigma_k.
\]
 Applying \eqref{eq:stable-interpolation-bound} to
$G(w)=G_1(\tau_k+w)$ gives
\[
 \left\|G_1\right\|_\infty=\left\|G_1\left(\tau_k+\cdot\right)\right\|_\infty\ge k
 \qquad(k\ge1),
\]
a contradiction. Therefore
\[
 W^2\notin AH^\infty\left(\mathbb C_+\right)+B_0H^\infty\left(\mathbb C_+\right).
\]

 Finally, transfer the construction back to the disk by setting
\[
 a=A\circ\Phi^{-1},\qquad b=B_0\circ\Phi^{-1},\qquad
 v=\tfrac14W\circ\Phi^{-1}.
\]
Equation \eqref{eqwabc} gives
$|v|\le|a|+|b|$. Composition with $\Phi^{-1}$ is an
isometric algebra isomorphism between the two $H^\infty$
spaces. Since $v^2\circ\Phi=W^2/16$, the nonmembership
just proved gives \eqref{eq:localized-nonmembership}. The lower bounds on
$K$ give \eqref{eq:localized-lower-bound} with $c=7/32$.
\end{proof}

\subsection{Intrinsic square obstruction}\label{sec:intrinsic-obstruction}

The localized example can now be reflected to obtain real-symmetric
data. Since $U$ is disjoint from its reflection, the lower bounds in
Lemma~\ref{lloctreil} preserve pointwise domination and allow us to
remove the reflected factor from an ideal relation by the corona theorem.
The slice regular extensions will then give Theorem~\ref{thm2-1}.

\begin{lem}
\label{lrealsymcounterexample}
There exist $P,Q,T\in H^\infty_{\mathbb R}(\mathbb D)$ such that
\[
 |T(z)|\le|P(z)|+|Q(z)|\le1,
 \qquad z\in\mathbb D,
\]
but
\[
 T^2\notin PH^\infty(\mathbb D)+QH^\infty(\mathbb D).
\]
\end{lem}

\begin{proof}
Take $a,b,v,c,U$ from Lemma~\ref{lloctreil}, with
$0<c\le1$. The set $U$ is disjoint from its reflection
$\left\{\bar z:z\in U\right\}$. Define
\[
 P=\frac{a\widehat a}{2},\qquad
 Q=\frac{b\widehat b}{2},\qquad
 T=\frac{cv\widehat v}{2}.
\]
These functions are real-symmetric, and $|P|+|Q|\le1$
because $\|a\|_\infty,\|b\|_\infty\le1$.
 To verify $|T|\le|P|+|Q|$, we use the lower bounds at whichever
of $z,\bar z$ lies outside $U$. If $z\in U$, then $\bar z\notin U$, so
$|\widehat a(z)|,|\widehat b(z)|\ge c$ and $|\widehat v(z)|\le1$.
It follows from  \eqref{eqlocdom} that
\[
 |P(z)|+|Q(z)|
 \ge\frac c2(|a(z)|+|b(z)|)
 \ge\frac c2|v(z)|\ge|T(z)|.
\]
If $z\notin U$, then $|a(z)|,|b(z)|\ge c$ and
$|v(z)|\le1$. Hence, applying \eqref{eqlocdom}
at $\bar z$ gives
\[
 |P(z)|+|Q(z)|
 \ge\frac c2\left(\left|\widehat a(z)\right|+\left|\widehat b(z)\right|\right)
 \ge\frac c2\left|\widehat v(z)\right|\ge|T(z)|.
\]

To prove the failure of ideal membership, argue by contradiction.
Suppose $T^2=PG_1+QG_2$ for some  $G_1,G_2\in H^\infty(\mathbb D)$.
Substituting the definitions of $P,Q,T$ gives
\[
 v^2\widehat v^{\,2}=aV_1+bV_2,\qquad
 V_1=\frac{2\widehat aG_1}{c^2},\quad
 V_2=\frac{2\widehat bG_2}{c^2}
 \in H^\infty(\mathbb D).
\]
 We use the corona theorem to remove the factor $\widehat v^{\,2}$
from this ideal relation.

The localized lower bounds in (\ref{eq:localized-lower-bound}) give
\[
 |a|+|b|+\left|\widehat v\right|^2
 \ge
 \begin{cases}
 c^2,&z\in U,\\
 2c\ge c^2,&z\notin U.
 \end{cases}
\]
By Carleson's corona theorem, there exist $u_1,u_2,u_0\in H^\infty(\mathbb D)$
such that $au_1+bu_2+\widehat v^{\,2}u_0=1$. Consequently,
\[
 \begin{aligned}
 v^2
 &=a\left(u_1v^2\right)+b\left(u_2v^2\right)
   +u_0\left(v^2\widehat v^{\,2}\right)\\
 &=a\left(u_1v^2+u_0V_1\right)+b\left(u_2v^2+u_0V_2\right),
 \end{aligned}
\]
contradicting \eqref{eq:localized-nonmembership}.
\end{proof}

\begin{proof}[Proof of Theorem \ref{thm2-1}]
Fix $I,J\in\mathbb S$
with $IJ=-JI$. Choose the real-symmetric functions $P,Q,T\in H^\infty_{\mathbb{R}}(\mathbb D)$ from
Lemma~\ref{lrealsymcounterexample}, and let
\[
 a=\operatorname{ext}(T),\qquad
 a_1=\operatorname{ext}(P),\qquad
 a_2=\operatorname{ext}(Q).
\]
By Lemma~\ref{lintrinreduc}, these extensions are
bounded and intrinsic, with constant modulus on each
conjugacy sphere. Consequently,
\begin{equation}\label{eqaqaa}
 |a(q)|\le\left|a_1(q)\right|+\left|a_2(q)\right|\le1,
 \qquad q\in\mathbb B.
\end{equation}
The restriction of $a^{\ast2}$ to the slice $\mathbb D_I$ is $T^2$,
with second splitting component zero. Lemma~\ref{lintrinreduc}
therefore gives
\[
 a^{\ast2}\in\mathcal I\left(a_1,a_2\right)
 \quad\Longleftrightarrow\quad
 T^2\in PH^\infty(\mathbb D)+QH^\infty(\mathbb D).
\]
Lemma~\ref{lrealsymcounterexample} excludes the right-hand
condition, which proves \eqref{eq:intrinsic-square-obstruction}.\end{proof}

\begin{proof}[Proof of Corollary~\ref{t2}]
For \textup{(i)}, restrict \eqref{et21} to any slice
$\mathbb D_I$. Wolff's ideal theorem gives
\[
 f_I^3\in\sum_{\ell=1}^n\left(f_\ell\right)_I H^\infty\left(\mathbb D_I\right).
\]
Since $f^{\ast3}$ is intrinsic and $\left(f^{\ast3}\right)_I=f_I^3$,
Lemma~\ref{lintrinreduc} gives \eqref{et22} with intrinsic
coefficients.

Part~\textup{(ii)} follows from Theorem \ref{thm2-1}, with zero generators appended when $n>2$.
\end{proof}

Rescaling the intrinsic square counterexample gives an obstruction
for the control function $h$ in the question posed in Section~\ref{s1}.
\begin{cor}\label{cor:quadratic-gauge-obstruction}
Let $n\ge2$ be an integer, and let
$h:[0,\infty)\to[0,\infty)$ be nondecreasing, with $h(0)=0$.
If
\[
 h(t)\ge ct^2,\qquad 0<t\le t_0,
\]
for some constants $c,t_0>0$, then there exist
$f,f_1,\ldots,f_n\in H_+^\infty(\mathbb B)$ such that
\[
 |f(q)|\le
 h\left(\sum_{\ell=1}^n\left|f_\ell(q)\right|\right),
 \qquad q\in\mathbb B,
\]
while
\[
 f\notin\mathcal I\left(f_1,\ldots,f_n\right).
\]
\end{cor}
\begin{proof}
Take $a,a_1,a_2$ from Theorem~\ref{thm2-1}, choose
$0<\lambda\le\min\left\{1,t_0\right\}$, and define
\[
 f=c\lambda^2a^{\ast2},\qquad
 f_1=\lambda a_1,\qquad f_2=\lambda a_2.
\]
Set $f_3=\cdots=f_n=0$ when $n>2$.
These functions belong to $H_+^\infty(\mathbb B)$.
 Here $S=\sum_{\ell=1}^n|f_\ell|=\lambda(|a_1|+|a_2|)\le t_0$.
Since $a$ is intrinsic, \eqref{eq:intrinsic-square-data} and the lower bound on $h$ give,
at points where $S>0$,
\[
 |f|=c\lambda^2|a|^2
 \le c\lambda^2\left(\left|a_1\right|+\left|a_2\right|\right)^2
 =cS^2\le h(S).
\]
At $S=0$, both $f$ and $h(S)$ vanish. Since $\lambda>0$,
\[
 \mathcal I\left(f_1,\ldots,f_n\right)=\mathcal I\left(a_1,a_2\right),\qquad
 f\in\mathcal I\left(f_1,\ldots,f_n\right)
 \ \Longrightarrow\
 a^{\ast2}=\frac{f}{c\lambda^2}\in\mathcal I\left(a_1,a_2\right),
\]
contrary to \eqref{eq:intrinsic-square-obstruction}.
\end{proof}

\medskip

\section{Sufficient conditions}
\label{sec:proof-gauge}

 We now prove the sufficient conditions in Theorem~\ref{thm2}.
Both parts use Treil's scalar ideal theorem. For general generators,
we apply it to the maximal minors of the matrix $M$ and then recover
a solution of $MX=H$ by the adjugate formula. For generators preserving
a common slice, we apply it directly to the two splitting components
of $f$.

\subsection{Scalar ideal estimates and normalization}

We use the following form of Treil's ideal theorem.

\begin{thm}[Treil {\cite[Theorem~0.1]{Tr07}}]\label{thmtdini}
Let $m\ge1$ and $A_1,\ldots,A_m,T\in H^\infty(\mathbb D)$.
Set
\[
 d(z)=\left(\sum_{\ell=1}^m\left|A_\ell(z)\right|^2\right)^{1/2},
\]
and assume that $d(z)\le1$ for every $z\in\mathbb D$.
Let $\eta:(0,1]\to[0,\infty)$ be bounded and nondecreasing,
with
\[
 \int_0^1\frac{\eta(t)}{t}\,\rd t<\infty,
\]
and set $\eta(0)=0$. If, for some constant $C\ge0$,
\[
 |T(z)|\le C d(z)^2\eta(d(z)),
 \qquad z\in\mathbb D,
\]
then
\[
 T\in\sum_{\ell=1}^m A_\ell H^\infty(\mathbb D).
\]
\end{thm}

This formulation follows by setting $\psi(x)=\eta(e^{-x/2})$
in Treil's theorem and rescaling $T$. Indeed,
\[
 \int_0^\infty\psi(x)\,\rd x
 =2\int_0^1\frac{\eta(t)}{t}\,\rd t<\infty.
\]

 Theorem~\ref{thmtdini} requires normalized generators and a bound on
the whole disk, whereas the decay assumption in Theorem~\ref{thm2}
is imposed only near $S=0$. The next lemma provides the normalization
and extends that local bound without changing the generated ideal.

\begin{lem}
\label{lgnormalization}
Let $\omega:(0,1]\to[0,\infty)$ be bounded and
nondecreasing, and suppose that
\[
 \int_0^1\frac{\omega(t)}{t}\,\rd t<\infty.
\]
Extend $\omega$ to $[0,\infty)$ by setting
$\omega(0)=0$ and $\omega(t)=\omega(1)$ for $t>1$.
Then the following statements hold.

\begin{enumerate}
\renewcommand{\labelenumi}{\textup{(\roman{enumi})}}

\item
For every $c,\alpha>0$, the function
$t\mapsto\omega\left(ct^\alpha\right)$ is bounded, nonnegative,
and nondecreasing on $(0,1]$, and satisfies
\[
 \int_0^1\frac{\omega\left(ct^\alpha\right)}{t}\,\rd t<\infty.
\]

\item
Let $n\ge1$ and
$f,f_1,\ldots,f_n\in H^\infty(\mathbb B)$, and put
\[
 S(q)=\sum_{\ell=1}^n\left|f_\ell(q)\right|,
 \qquad q\in\mathbb B.
\]
Suppose that $h:[0,\infty)\to[0,\infty)$ is
nondecreasing, $h(0)=0$, and
\[
 |f(q)|\le h(S(q)),
 \qquad q\in\mathbb B.
\]
Assume also that there exist constants
$C_0>0$, $r>0$, and $t_0\in(0,1]$ such that
\[
 h(t)\le C_0t^r\omega(t),
 \qquad 0<t\le t_0.
\]

Define
\[
 L=\max\left\{1,\sum_{\ell=1}^n
                       \left\|f_\ell\right\|_\infty\right\}.
\]
Then there exist a constant $C_1>0$ and a bounded,
nonnegative, nondecreasing function
$\omega_1:(0,1]\to[0,\infty)$ satisfying
\[
 \int_0^1\frac{\omega_1(t)}{t}\,\rd t<\infty
\]
such that, after setting $\omega_1(0)=0$,
\begin{equation}\label{eq:normalized-gauge}
 \sum_{\ell=1}^n\left|\frac{f_\ell(q)}{L}\right|
 =\frac{S(q)}{L}\le1,
 \qquad
 |f(q)|\le
 C_1\left(\frac{S(q)}{L}\right)^r
 \omega_1\left(\frac{S(q)}{L}\right),
 \qquad q\in\mathbb B.
\end{equation}
Moreover,
\[
 \mathcal I\left(\frac{f_1}{L},\ldots,\frac{f_n}{L}\right)
 =\mathcal I\left(f_1,\ldots,f_n\right).
\]
\end{enumerate}

Statement \textup{(ii)} also holds for
$f,f_1,\ldots,f_n\in H^\infty(\mathbb D)$, with
$q\in\mathbb B$ replaced by $z\in\mathbb D$ and
slice products replaced by pointwise products.
\end{lem}

\begin{proof}
\textup{(i)} Monotonicity and the Dini condition of $\omega$ imply
$\lim_{t\downarrow0}\omega(t)=0$.
The boundedness and monotonicity of the function
$t\mapsto\omega\left(ct^\alpha\right)$ follow from
the same properties of $\omega$. The change of variables
$s=ct^\alpha$ gives
\[
 \int_0^1\frac{\omega\left(ct^\alpha\right)}t\,\rd t
 =\frac1\alpha\int_0^c\frac{\omega(s)}s\,\rd s<\infty.
\]
The integral is finite near zero by assumption and on any
remaining interval by boundedness of the extended function.

 \textup{(ii)}   Note that $0\le S(q)\le L$.
The decay hypothesis controls $f$ where $S\le t_0$; elsewhere,
the bound $h(S)\le h(L)$ suffices. Accordingly, define
\[
 \omega_1(t)=
 \begin{cases}
 \omega(Lt),&0<t\le t_0/L,\\
 \max\{1,\omega(1)\},&t_0/L<t\le1,
 \end{cases}
 \qquad \omega_1(0)=0.
\]
This function is bounded and nondecreasing, since the constant
in the second branch is at least $\omega\left(t_0\right)$. Moreover,
\[
 \int_0^1\frac{\omega_1(t)}t\,\rd t
 =\int_0^{t_0}\frac{\omega(s)}s\,\rd s
  +\max\{1,\omega(1)\}\log\frac{L}{t_0}<\infty.
\]
For $t=S(q)/L\in[0,1]$, if
$0<t\le t_0/L$, then
\[
 |f(q)|\le h(Lt)
 \le C_0L^r t^r\omega(Lt)
 =C_0L^r t^r\omega_1(t).
\]
If $t_0/L<t\le1$, then $\omega_1(t)\ge1$ and
\[
 |f(q)|\le h(Lt)\le h(L)
 \le h(L)\left(L/t_0\right)^r t^r\omega_1(t).
\]
At $t=0$, the assumption $h(0)=0$ gives $f(q)=0$.
Thus \eqref{eq:normalized-gauge} holds with
\[
 C_1=\max\left\{C_0L^r,\,h(L)\left(L/t_0\right)^r\right\}>0.
\]
Finally,
$\left(f_\ell/L\right)\ast g_\ell=f_\ell\ast\left(g_\ell/L\right)$, and real
scaling is a bijection of $H^\infty(\mathbb B)$, so the ideals
coincide. The same proof applies to complex functions.
\end{proof}

\subsection{Maximal minors and a two-row estimate}

For $M\in\mathbb C^{2\times r}$, $r\ge2$, write
$0\le a\le b=\|M\|_{\mathrm{op}}$ for its singular values.
Its maximal minors are the determinants of the $2\times2$
submatrices formed by choosing two columns.

The Cauchy--Binet formula gives the following identity; see \cite[Section~0.8.7]{HJ}.
\begin{lem}\label{lcbinet}
Let $M\in\mathbb C^{2\times r}$, where $r\ge2$.
For $1\le\alpha<\beta\le r$, let $ M^{[\alpha,\beta] }$ be
the $2\times2$ matrix formed from columns $\alpha,\beta$
of $M$. Then
\begin{equation}\label{eqcbinet}
 \det\left(MM^*\right)
 =\sum_{1\le\alpha<\beta\le r}\left|\det M^{[\alpha,\beta]}\right|^2.
\end{equation}
Consequently, if $0\le a\le b$ are the singular values
of $M$, then
\[
 D=\left(\sum_{\alpha<\beta}\left|\det M^{[\alpha,\beta]}\right|^2\right)^{1/2}
 =ab.
\]
\end{lem}

 To pass from the minors to a matrix solution, we use the classical
adjugate identity
\[
 A\operatorname{adj}(A)=(\det A)I_2,\qquad A\in\mathbb C^{2\times2}.
\]
The entries of $\operatorname{adj}(A)$ are linear in those of $A$;
thus this construction preserves boundedness and holomorphicity.
Compare~\cite{HT}.

\begin{lem}
\label{lemmlifting}
Let $M$ be a bounded holomorphic $2\times r$ matrix on
$\mathbb D$, where $r\ge2$. For $1\le\alpha<\beta\le r$,
let $E_{\alpha\beta}$ be the $r\times2$ matrix whose entries
are zero except for a one in positions $(\alpha,1)$ and
$(\beta,2)$. Thus $ M^{[\alpha,\beta]}=ME_{\alpha\beta }$ consists of
columns $\alpha,\beta$ of $M$. Set
\[
 {d_{\alpha\beta}=\det M^{[\alpha,\beta]}.}
\]
Let $H$ be a bounded holomorphic $2\times2$ matrix on
$\mathbb D$. If
\[
 H_{\mu\nu}\in\sum_{1\le\alpha<\beta\le r}
                  d_{\alpha\beta}H^\infty(\mathbb D),
 \qquad 1\le\mu,\nu\le2,
\]
then there exists a bounded holomorphic $r\times2$ matrix
$X$ on $\mathbb D$ such that $MX=H$.
For the block matrices in \eqref{eqblomat}, this implies
$f\in\mathcal I\left(f_1,\ldots,f_n\right)$.
\end{lem}

\begin{proof}
By assumption, there exist
$k_{\alpha\beta}^{\mu\nu}\in H^\infty(\mathbb D)$ such that
\[
 H_{\mu\nu}=\sum_{\alpha<\beta}
                d_{\alpha\beta}k_{\alpha\beta}^{\mu\nu}.
\]
Set $K_{\alpha\beta}=\left(k_{\alpha\beta}^{\mu\nu}\right)_{\mu,\nu=1}^2$.
Then $H=\sum_{\alpha<\beta}d_{\alpha\beta}K_{\alpha\beta}$.
 The adjugate identity expresses each term
$d_{\alpha\beta}K_{\alpha\beta}$ as a matrix product with $M$ on the left.
Define
\begin{equation}\label{eXeadj}
  X=\sum_{\alpha<\beta}E_{\alpha\beta}
       \operatorname{adj}\left(M^{[\alpha,\beta]}\right)K_{\alpha\beta}.
\end{equation}
The entries of $X$ are bounded and holomorphic, and
\[
\begin{aligned}
 MX
 &=\sum_{\alpha<\beta}
      M^{[\alpha,\beta]}\operatorname{adj}\left(M^{[\alpha,\beta]}
      \right)K_{\alpha\beta} =\sum_{\alpha<\beta}d_{\alpha\beta}K_{\alpha\beta}=H.
\end{aligned}
\]
For the matrices in \eqref{eqblomat}, write
\[
 X=\begin{pmatrix}X_1\\\vdots\\X_n\end{pmatrix},
 \qquad X_\ell\in H^\infty(\mathbb D)^{2\times2}.
\]
Then $MX=H$ becomes
$\sum_{\ell=1}^n\chi\left(f_\ell\right)X_\ell=\chi(f)$, so
Lemma~\ref{lmsymmetry} gives $f\in\mathcal I\left(f_1,\ldots,f_n\right)$.
\end{proof}

 By Lemma~\ref{lemmlifting}, to construct a bounded holomorphic solution
of $MX=H$, it suffices to show that every entry of $H$ belongs to the
ideal generated by the maximal minors of $M$.  To apply Theorem~\ref{thmtdini} to the minors,
we need a bound in terms of their Euclidean norm $D$.
The next lemma obtains it from estimates for $uH$ and $uM$.

\begin{lem} \label{ltrowestimate}
Let $M\in\mathbb C^{2\times r}$, where $r\ge2$, satisfy
$\|M\|_{\mathrm{op}}\le1$, with singular values
$0\le a\le b\le1$, and put $D=ab$.
Let $H\in\mathbb C^{2\times2}$, and let
$\omega:(0,1]\to[0,\infty)$ be bounded and nondecreasing.
Extend $\omega$ to $[0,\infty)$ by setting $\omega(0)=0$
and $\omega(t)=\omega(1)$ for $t>1$.
Suppose that there are constants $K,L>0$ such that, for
every unit row $u\in\mathbb C^2$,
\begin{equation}\label{etworow}
 \|uH\|\le K\|uM\|^4\omega(L\|uM\|).
\end{equation}
Then
\begin{equation}\label{eq:two-row-conclusion}
 \|H\|_{\mathrm{op}}
 \le8Kba^3\omega(\sqrt2La)
 \le8KD^2\omega\left(\sqrt2L\sqrt D\right).
\end{equation}
In particular, $H=0$ whenever $\operatorname{rank}M<2$.
\end{lem}

\smallskip
\begin{proof}
Choose orthonormal rows $e_1,e_2$ with
$e_1MM^*=a^2e_1$ and $e_2MM^*=b^2e_2$. Then
\[
 \left\|e_1M\right\|=a,\qquad \left\|e_2M\right\|=b,\qquad
 \left(e_1M\right)\left(e_2M\right)^*=0.
\]
 Suppose $a>0$.  To estimate $e_2H$, choose
$t=a/b$ so that the contributions of $e_1M$ and $te_2M$ have equal norms.
Set $u_\pm=(e_1\pm te_2)/\sqrt{1+t^2}$. Orthogonality gives
\[
 \left\|u_\pm\right\|=1,\qquad
 \left\|u_\pm M\right\|^2=\frac{a^2+t^2b^2}{1+t^2}
 =\frac{2a^2}{1+t^2}\le2a^2.
\]
By \eqref{etworow} and monotonicity, we have
\[
 \left\|e_1H\right\|\le Ka^4\omega(La),\qquad
 \left\|u_\pm H\right\|\le4Ka^4\omega(\sqrt2La).
\]
Since $e_2=\sqrt{1+t^2}\left(u_+-u_-\right)/(2t)$, it follows that
\[
 \left\|e_2H\right\|
 \le\frac{\sqrt{1+t^2}}{2t}
        \left(\left\|u_+H\right\|+\left\|u_-H\right\|\right) \le4\sqrt2Kba^3\omega(\sqrt2La).
\]
Every unit row is $u=c_1e_1+c_2e_2$, with
$\left|c_1\right|^2+\left|c_2\right|^2=1$. Consequently,
\[
 \begin{aligned}
 \|H\|_{\mathrm{op}}
 &\le\left\|e_1H\right\|+\left\|e_2H\right\|
 \le K\left(a^4+4\sqrt2ba^3\right)\omega(\sqrt2La)\\
 &\le8Kba^3\omega(\sqrt2La)
 \le8KD^2\omega\left(\sqrt2L\sqrt D\right),
 \end{aligned}
\]
where $a^4\le ba^3$, $1+4\sqrt2<8$,
$a\le\sqrt{ab}=\sqrt D$, and $ba^3\le a^2b^2=D^2$.

 When $a=0<b$, the same idea uses rows approaching $e_1$.
Equation~\eqref{etworow} first gives $e_1H=0$.
For $u_t=(e_1+te_2)/\sqrt{1+t^2}$, $0<t\le1$, it then gives
\[
 \frac{t\left\|e_2H\right\|}{\sqrt{1+t^2}}
 \le K\frac{t^4b^4}{\left(1+t^2\right)^2}
          \omega\left(\frac{Ltb}{\sqrt{1+t^2}}\right).
\]
Hence
\[
 \left\|e_2H\right\|
 \le\frac{Kb^4t^3}{\left(1+t^2\right)^{3/2}}
                 \sup_{s\ge0}\omega(s)
 \longrightarrow0\qquad(t\downarrow0),
\]
so $H=0$. If $b=0$, \eqref{etworow} directly gives $uH=0$
for every unit row $u$, and again $H=0$.
\end{proof}

\begin{rem}
The exponent $4$ in Lemma~\ref{ltrowestimate} reflects the estimate
needed to apply the scalar ideal theorem to the maximal minors.
For $0<t<1$ and $p>0$, let
\[
 M_t=\left[\,tI_2\ \ 0_{2\times2}\,\right],\qquad H_t=t^pI_2.
\]
For every unit row $u$, we have $\|uH_t\|=\|uM_t\|^p=t^p$,
whereas $D_t=\sqrt{\det(M_tM_t^*)}=t^2$. Thus
\[
 \frac{\|H_t\|_{\mathrm{op}}}{D_t^2}
 =t^{p-4}\longrightarrow\infty\qquad(t\downarrow0,\ p<4).
\]
Consequently, row domination with exponent $p<4$ does not imply
a uniform estimate $\|H_t\|_{\mathrm{op}}\le C D_t^2\eta(D_t)$
with bounded $\eta$. This concerns the minor estimate: the corresponding
matrix equation itself has the explicit solution
\[
 M_t\begin{pmatrix}t^{p-1}I_2\\0_{2\times2}\end{pmatrix}=H_t.
\]
\end{rem}

\subsection{Proof of the sufficient conditions}

\begin{proof}[Proof of Theorem~\ref{thm2}]
 \textup{(i)}  Apply Lemma~\ref{lgnormalization}\textup{(ii)} with $r=4$.
Retain the notation $f_\ell,S,\omega$ after normalization.
The ideal is unchanged, and
\begin{equation}\label{eqsqsf}
 S(q)=\sum_{\ell=1}^n\left|f_\ell(q)\right|\le1,\qquad
 |f(q)|\le C_0S(q)^4\omega(S(q)),\qquad q\in\mathbb B,
\end{equation}
for some $C_0>0$. Extend $\omega$ to $[0,\infty)$ as in
Lemma~\ref{lgnormalization}, in particular with $\omega(0)=0$.

Fix $I,J\in\mathbb S$ with $IJ=-JI$ and form $H,M$ as in
\eqref{eqblomat}. For a unit row $u=\left(u_1,u_2\right)$, put
$s=u_1+u_2J$ and $q=s^{-1}zs$. Equation \eqref{eqbrnorm}
and Cauchy--Schwarz give
\[
 \|uH(z)\|=|f(q)|,\qquad
 \|uM(z)\|\le S(q)\le\sqrt n\,\|uM(z)\|.
\]
Since $S\le1$, we obtain $\|M(z)\|_{\mathrm{op}}\le1$ and
\begin{equation}\label{eq:proof-gauge-row}
 \|uH(z)\|
 \le C_0S(q)^4\omega(S(q)) \le C_0n^2\|uM(z)\|^4
           \omega\left(\sqrt n\,\|uM(z)\|\right).
\end{equation}
 We now convert this row estimate into the scalar estimate needed
for the minors in Lemma~\ref{lemmlifting}.
 For $1\le\alpha<\beta\le2n$, put
\[
M^{[\alpha,\beta]}=ME_{\alpha\beta},\qquad
 d_{\alpha\beta}=\det M^{[\alpha,\beta]},\qquad
 D(z)=\left(\sum_{\alpha<
 \beta}\left|d_{\alpha\beta}(z)\right|^2\right)^{1/2},
\]
where $E_{\alpha\beta}$ is the column-selection matrix from
Lemma~\ref{lemmlifting}. The minors $d_{\alpha\beta}$ are bounded
and holomorphic. By Lemma~\ref{lcbinet},
\[
 D(z)=a(z)b(z)\le1,
\]
where $0\le a(z)\le b(z)\le1$ are the singular values of $M(z)$.
Applying Lemma~\ref{ltrowestimate} to
\eqref{eq:proof-gauge-row}, with $C_0n^2$ and $\sqrt n$ as
the two constants in \eqref{etworow}, gives
\begin{equation}\label{epgmibound}
 \|H(z)\|_{\mathrm{op}}
 \le8C_0n^2D(z)^2\omega\left(\sqrt{2n}\sqrt{D(z)}\right).
\end{equation}
 The bound \eqref{epgmibound} has the form required by
Theorem~\ref{thmtdini} once we check the Dini condition for
$\eta(t)=\omega(\sqrt{2n}\sqrt t)$, $0\le t\le1$.
This function is bounded, nonnegative and nondecreasing, with $\eta(0)=0$, and
\begin{equation}\label{eqdini}
 \int_0^1\frac{\eta(t)}t\,\rd t
 =2\int_0^{\sqrt{2n}}\frac{\omega(s)}s\,\rd s =2\int_0^1\frac{\omega(s)}s\,\rd s
    +\omega(1)\log(2n)<\infty.
\end{equation}
For each $1\le\mu,\nu\le2$,
\[
 \left|H_{\mu\nu}(z)\right|\le\|H(z)\|_{\mathrm{op}}
 \le8C_0n^2D(z)^2\eta(D(z)).
\]
Theorem~\ref{thmtdini}, applied to the scalar tuple
$\left(d_{\alpha\beta}\right)_{\alpha<\beta}$, therefore gives
\[
 H_{\mu\nu}\in\sum_{\alpha<\beta}
                    d_{\alpha\beta}H^\infty(\mathbb D),
 \qquad1\le\mu,\nu\le2.
\]
 By Lemma~\ref{lemmlifting}, these four scalar representations give
a bounded holomorphic solution of $MX=H$. Symmetrizing its $2\times2$
blocks as in Lemma~\ref{lmsymmetry} then gives
$f\in\mathcal I(f_1,\ldots,f_n)$.

  \textup{(ii)}  Lemma~\ref{lintrinreduc}\textup{(ii)} reduces the ideal
equation to two scalar equations on the common slice.
Use Lemma~\ref{lgnormalization} with $r=2$ to arrange
\[
 S\le1,\qquad |f|\le C_0S^2\omega(S).
\]
Extend $\omega$ constantly beyond one and set $\omega(0)=0$,
as in part~\textup{(i)}. Choose a common slice $\mathbb D_I$
and an orthogonal unit $J$. Identify $\mathbb D_I$ with $\mathbb D$, keep the notation
$f_\ell$ for the restrictions, and write $f_I=U+VJ$.
Set
\[
 d(z)=\left(\sum_{\ell=1}^n\left|f_\ell(z)\right|^2\right)^{1/2},
 \qquad d\le S\le\min\left\{1,\sqrt n\,d\right\}.
\]
By \eqref{eq:splitting-bounds}, we have
\begin{equation}\label{eq:proof-common-slice-bound}
 \max\{|U(z)|,|V(z)|\}
 \le\left|f_I(z)\right|\le C_0S(z)^2\omega(S(z))
 \le C_0n d(z)^2\omega\left(\sqrt n\,d(z)\right).
\end{equation}
 Moreover,
\[
 \int_0^1\frac{\omega\left(\sqrt n\,t\right)}t\,\rd t
 =\int_0^1\frac{\omega(s)}s\,\rd s
    +\tfrac12\omega(1)\log n<\infty.
\]
Applying Theorem~\ref{thmtdini} to $U$ and $V$, then
Lemma~\ref{lintrinreduc}, gives
\[
 U,V\in\sum_{\ell=1}^n f_\ell H^\infty(\mathbb D)
 \quad\Longrightarrow\quad f\in\mathcal I\left(f_1,\ldots,f_n\right).
\]
 For intrinsic data, Lemma~\ref{lintrinreduc}\textup{(iii)} gives intrinsic
coefficients in both \textup{(i)} and \textup{(ii)}.
\end{proof}

\begin{proof}[Proof of Corollary~\ref{thm:power-ideal}]
Set $h(t)=Ct^p$. For arbitrary generators and $p>4$, use
$\omega(t)=t^{p-4}$ in Theorem~\ref{thm2}\textup{(i)};
for a common slice and $p>2$, use $\omega(t)=t^{p-2}$ in
part~\textup{(ii)}. These functions are bounded and
nondecreasing on $(0,1]$, and
\[
 \int_0^1\frac{t^{p-r}}t\,\rd t=\frac1{p-r}<\infty
 \qquad(p>r,\ r\in\{2,4\}).
\]
The same theorem supplies intrinsic coefficients for intrinsic data.

For the negative assertion, let $n\ge2$ and take $a,a_1,a_2$ from
Theorem \ref{thm2-1}, put $f=a^{\ast2}$,
$f_1=a_1$, $f_2=a_2$, and append zero generators if $n>2$.
Then $S=\left|a_1\right|+\left|a_2\right|\le1$ and
\[
 |f|=|a|^2\le S^2\le S^p\qquad(1\le p\le2),
 \qquad
 f\notin\mathcal I\left(a_1,a_2\right)=\mathcal I\left(f_1,\ldots,f_n\right).
\]
The same data work for all $p\in[1,2]$.
\end{proof}

\medskip

\subsection*{Funding} Xu is partially supported by  China Scholarship Council (No. 202506690052).
 Yang is partially supported by   Guangdong Basic and Applied Basic Research Foundation (No. 2023A1515110597).
 Yuan  is partially supported by   Guangdong Basic and Applied Basic Research Foundation (No. 2026A1515012908).

\section*{Statements and Declarations}

\subsection*{Data availability}
No datasets were generated or analyzed during this study.

\subsection*{Conflict of interest}
The authors declare that they have no conflicts of interest.

\subsection*{Use of AI tools}

The mathematical calculations and derivations in this paper were carried out by the authors. After completing the initial draft, the authors used OpenAI's ChatGPT (GPT-6 Astra) and DeepSeek (DeepSeek-V4.1) for English-language editing and to review the calculations and derivations. The authors directed these tasks, reviewed the tools' suggestions, and independently verified the material retained in the manuscript. All authors take responsibility for the final content.

\end{document}